\documentclass[11pt, twosided]{article}
\usepackage{amsmath,amscd,amsfonts,eucal,latexsym,amssymb,mathrsfs}
\usepackage{amsxtra, amsthm, calc}
\usepackage[usenames]{color}
\theoremstyle{definition}
\newtheorem{Definition}{Definition}[section]
\theoremstyle{plain}
\newtheorem{Theorem}[Definition]{Theorem}
\newtheorem{Proposition}[Definition]{Proposition}
\newtheorem{Lemma}[Definition]{Lemma}
\newtheorem{Corollary}[Definition]{Corollary}
\theoremstyle{remark}
\newtheorem{Remark}[Definition]{Remark}
\newtheorem{Remarks}[Definition]{Remarks}
\newtheorem{example}[Definition]{Example}
\newcounter{remcount}

\newenvironment{remlist}{\begin{list}{(\roman{remcount})}{\usecounter{remcount}\setlength{\leftmargin}{0 cm}\setlength{\rightmargin}{0 cm}\setlength{\topsep}{0 cm}\setlength{\itemsep}{0 cm}\setlength{\parsep}{\parskip}\setlength{\labelwidth}{1 cm}\setlength{\labelsep}{0.5 em}\setlength{\itemindent}{1 cm + 0.5 em}}}{\end{list}}
\definecolor{Red}{rgb}{1,0,0}
\definecolor{Green}{rgb}{0,0.65,0}
\def\cA{{\cal A}}
\def\cB{{\cal B}}
\def\cC{{\cal C}}
\def\cD{{\cal D}}

\def\cF{{\cal F}}
\def\cG{{\cal G}}
\def\cH{{\cal H}}
\def\cI{{\cal I}}

\def\cN{{\cal N}}

\def\bC{{\mathbb C}}

\def\bN{{\mathbb N}}

\def\bR{{\mathbb R}}

\def\bZ{{\mathbb Z}}

\def\a{\alpha}
\def\b{\beta}
\def\g{\gamma}        
\def\d{\delta}        \def\D{\Delta}

\def\k{\kappa}
\def\l{\lambda}       \def\L{\Lambda}
\def\m{\mu}

\def\s{\sigma}
    
\def\o{\omega}        \def\O{\Omega}

\def\fC{{\mathfrak C}}

\def\fG{{\mathfrak G}}

\def\supp{{\text{supp}}}
\def\ad{{\text{Ad}}}
\newcommand{\rest}{\upharpoonright}
\newcommand{\id}{1}
\newcommand{\uF}{\underline{F}}

\newcommand{\uA}{\underline{A}}
\newcommand{\uB}{\underline{B}}
\newcommand{\uFl}{\underline{F}_{\lambda}}

\newcommand{\uFF}{\underline{\mathfrak F}}
\newcommand{\uAA}{\underline{\mathfrak A}}
\newcommand{\uBB}{\underline{\mathfrak B}}
\newcommand{\ua}{\underline{\alpha}}

\newcommand{\uE}{\underline{E}}
\newcommand{\uooi}{\underline{\o}_{0,\iota}}

\newcommand{\Foi}{{\mathcal{F}_{0,\iota}}}
\newcommand{\Aoi}{\mathcal{A}_{0,\iota}}
\newcommand{\Boi}{\mathcal{B}_{0,\iota}}
\newcommand{\uo}{\underline{\omega}}

\newcommand{\ub}{\smash{\underline{\beta}}}
\newcommand{\aoi}{\alpha^{(0, \iota)}}
\newcommand{\boi}{\beta^{(0,\iota)}}
\newcommand{\poi}{\pi_{0,\iota}}
\newcommand{\Hoi}{{\cal H}_{0,\iota}}
\newcommand{\Ooi}{\Omega_{0,\iota}}
\newcommand{\ooi}{\omega_{0,\iota}}

\newcommand{\Uoi}{U_{0,\iota}}
\newcommand{\Voi}{V_{0,\iota}}

\newcommand{\Goi}{G_{0,\iota}}

\newcommand{\Eoi}{E_{0,\iota}}
\newcommand{\eoi}{e_{0,\iota}}
\newcommand{\cFm}{\cF^{(m)}}
\newcommand{\cAm}{\cA^{(m)}}
\newcommand{\am}{\a^{(m)}}
\newcommand{\om}{\o^{(m)}}
\newcommand{\cFz}{\cF^{(0)}}
\newcommand{\cAz}{\cA^{(0)}}
\newcommand{\az}{\a^{(0)}}
\newcommand{\oz}{\o^{(0)}}

\newcommand{\cFmoir}{\cFm_{0,\iota;r}}
\newcommand{\cAmoi}{\cAm_{0,\iota}}
\newcommand{\cAmoir}{\cAm_{0,\iota;r}}
\newcommand{\amoi}{\a^{(m;\,0,\iota)}}
\newcommand{\omoi}{\om_{0,\iota}}
\newcommand{\uAm}{\uAA^{(m)}}
\newcommand{\uFm}{\uFF^{(m)}}
\newcommand{\bmu}{\boldsymbol{\mu}}
\newcommand{\cFbm}{\cF^{(\bmu)}}
\newcommand{\cAbm}{\cA^{(\bmu)}}
\newcommand{\cFbz}{\cF^{(\boldsymbol{0})}}
\newcommand{\cAbz}{\cA^{(\boldsymbol{0})}}
\newcommand{\cFbmoi}{\cF^{(\bmu)}_{0,\iota}}
\newcommand{\cAbmoi}{\cA^{(\bmu)}_{0,\iota}}

\newcommand{\Ebmoi}{E^{(\bmu)}_{0,\iota}}

\newcommand{\pbmoi}{\pi^{(\bmu)}_{0,\iota}}
\newcommand{\cHbz}{\cH^{(\boldsymbol{0})}}
\newsymbol\bt 1202  

\newcommand{\norm}[1]{\left\lVert #1 \right\rVert}

\newcommand{\hide}[1]{} 

\newcommand{\cont}{\subseteq}

\newcommand{\Id}{{\bf 1}}

\newcounter{propcount}

\newlength{\maxlabelwidth}
\newenvironment{proplist}[2][1]{\begin{list}{(\roman{propcount})}{\usecounter{propcount}\setcounter{propcount}{#2}\settowidth{\maxlabelwidth}{\textit{(\roman{propcount})}}\setcounter{propcount}{#1-1}\setlength{\leftmargin}{\maxlabelwidth+
0.5 em}\setlength{\rightmargin}{0 cm}\setlength{\topsep}{0
cm}\setlength{\itemsep}{0
cm}\setlength{\parsep}{\parskip}\setlength{\labelwidth}{\maxlabelwidth}\setlength{\labelsep}{0.5
em}\setlength{\itemindent}{0 cm}}}{\end{list}}

\newcommand{\inst}[1]{$^\textrm{#1}$ }
\newcommand{\email}[1]{e-mail: #1}

\begin{document}
\title{Scaling limit for subsystems and \\
Doplicher-Roberts reconstruction}
\author{Roberto Conti\inst{1}
\and Gerardo Morsella\inst{2}}

\date{
\parbox[t]{0.9\textwidth}{\footnotesize{%
\begin{enumerate}
\renewcommand{\theenumi}{\arabic{enumi}}
\renewcommand{\labelenumi}{\theenumi}
\item Mathematics, School of Mathematical and Physical Sciences, The University of Newcastle, NSW 2308, Australia, \email{Roberto.Conti@newcastle.edu.au}
\item Scuola Normale Superiore di Pisa, Piazza dei Cavalieri, 7, 56126 Pisa, Italy, \email{gerardo.morsella@sns.it}
\end{enumerate}
}}
\\
\vspace{\baselineskip}
\today}

\maketitle

\begin{abstract}
Given an inclusion $\cB \subset \cF$ of (graded) local nets, we analyse the structure of the corresponding inclusion of scaling limit nets $\cB_0 \subset \cF_0$, giving conditions, fulfilled in free field theory, under which the unicity of the scaling limit of $\cF$ implies that of the scaling limit of $\cB$. As a byproduct, we compute explicitly the (unique) scaling limit of the fixpoint nets of scalar free field theories. In the particular case of an inclusion $\cA \subset \cB$ of local nets with the same canonical field net $\cF$, we find sufficient conditions which entail the equality of the canonical field nets of $\cA_0$ and $\cB_0$.
\end{abstract}

\section{Introduction}

Local quantum physics is an approach to Quantum Field Theory (QFT) based only on observable quantities~\cite{Haag}.
It has been very successful in the  mathematical description of superselection
sectors and of the global gauge group of a given QFT~\cite{DR}. Also, the mathematical
tools that are available in this setting are well suited for providing a
detailed analysis of subsystems, an issue that is central in order to obtain an
intrinsic description of the observable system one starts with~\cite{CDR, CC1, CC2}.
In another direction, the recently proposed algebraic approach to the renormalization group~\cite{BV} (see also section~\ref{sec:subsystem}) has
opened the possibility of studying the short distance limit in the local quantum physics
framework, and has started
to convey new insight into our
understanding of physically relevant issues such as confinement of colour charges
and renormalization of pointlike fields~\cite{Buc1, DMV, BDM1}.

\medskip

Dealing with the general problem of understanding the scaling limit $\cA_0$ of a given local net $\cA$,
it is natural to ask whether there exists an efficient way to compute it in practical situations.

Loosely speaking, starting with a given local net, one would like to mod out the degrees of freedom that play no role at short scale, and obtain a smaller and hopefully simpler net which has the scaling limit of the net we started with.
In turn, it is unlikely that a local net always contains a ``large'' subnet containing the whole information about scaling,
however
the notion of convergent scaling limit that we use at some crucial point of the main text is clearly an evolution of this naive idea.
In fact, the very concept of the scaling algebra involves some redundancy
in the choice of the scaling functions (as shown by the presence of
a big kernel of the scaling limit representation), so that one could
expect that, at least in particular cases, the consideration of some appropriate subalgebra
of the scaling algebra would suffice.

It might also be the case that one knows that the given local net can be realized as fixpoints of a larger net, and then wonder if the scaling limit of the fixpoints can be computed as the fixpoint of the scaling limit.

In both cases, we are thus led to the problem of comparing the scaling limit of a system with that of a subsystem,
and this paper came out as an attempt to understand this relationship.

\medskip
It has been shown in \cite{BV} that, for a given theory $\cA$, there are only three possibilities:
either $\cA$ has a trivial scaling limit,
or a unique non trivial one,
or several non-isomorphic ones.
In the case of a subsystem $\cA \subset \cB$ the situation is slightly more complicated,
but one of our most basic observations is that there always exists a bijective correspondence
between the sets of scaling limit states of $\cA$ and $\cB$, and that for the
corresponding scaling limit nets $\cA_0$ and $\cB_0$ one has a subsystem
$\cA_0 \subset \cB_0$. In this situation it seems natural to expect that if $\cB$ has a unique scaling limit,
the same should be true for $\cA$, at least under suitable assumptions. In section~\ref{sec:subsystem}
we provide a criterion for this to happen, and use it to show that fixed point nets
in free field theory have a unique scaling limit.

Another aspect of the problem is to study the inclusions $\cA_0 \subset \cB_0$, and for instance one can ask whether it is possible to find necessary/sufficient conditions on $\cA$ and $\cB$ ensuring that
\begin{equation}\label{eq:equal}
\cA_0 = \cB_0
\end{equation}
for every scaling limit state.

At first sight,
one could expect that the situation becomes somehow easier to handle if one knows
that
($\cA(O)$ and $\cB(O)$ are factors and)
 $[\cB:\cA] < \infty$, but it has not really become important to employ this
condition yet.
In turn,
the index of an inclusion is not necessarily preserved in the scaling limit,
but there are cases in which the inequality $[\cB_0:\cA_0] \leq [\cB:\cA]$ holds true.
For instance, consider the free massive scalar field and its $\mathbb{Z}_2$-fixpoints.
In this case, after scaling the index remains the same, as shown in section~\ref{sec:subsystem}.
On the other hand, tensorizing  the Lutz model~\cite{L} with a massive free field
one gets that the index of the scaling limit is smaller than the original one~\cite{DM}.
Also, the relation $[\cB_0:\cA_0] \leq [\cB:\cA]$ is compatible with equation~\eqref{eq:scalingfixpoint}
below.

Tensor products provide simple examples of subsystems, for which
some questions can be answered.
For instance, let us assume that $\cA_2$ has trivial scaling limit.\footnote{We remark that the role played by
Haag duality in relation to the triviality of the scaling limit is not completely understood, as no examples of nets satisfying it and having trivial scaling limit are known to date.} It is then natural to ask
under which conditions the scaling limits of $\cB := \cA_1 \otimes \cA_2$ and $\cA := \cA_1 \otimes {\mathbb C} \simeq \cA_1$ satisfy equation~\eqref{eq:equal}. A set of sufficient conditions for this to happen, expressed in terms of
nuclearity properties, has been found in~\cite{DM}. The fact that nuclearity plays a role in this context
is not surprising, as it provides invariants which depend on the localization region, and therefore
should be able to encode the fact that $\cA(O)$ and $\cB(O)$ become ``closer'' at small distances.

Notice also that a certain (graded) tensor product decomposition
plays a critical role in the classification of subsystems in \cite{CC1,CC2}.

This work heavily relies on
the DR-reconstruction \cite{DR}: given an observable net $\cA$, there exist a canonical field net
$\cF(\cA)$ and a compact group $G(\cA)$ of automorphisms of $\cF(\cA)$ such that $\cA=\cF(\cA)^{G(\cA)}$ (therefore $\cA$ is a subsystem of $\cF(\cA)$).
Some functoriality aspects of the reconstruction have been investigated in \cite{CDR},
and a classification result for subsystems of $\cF(\cA)$ has been obtained in \cite{CC1,CC2}.

The study of the scaling limits of $\cA$ and $\cF(\cA)$ is discussed in \cite{DMV}.
In typical cases, it holds
\begin{align}
\cF(\cA)_0 &= \cF(\cA_0)^H,\label{eq:scalingfieldfixpoint} \\
\cA_0 &= \cF(\cA)_0^{G(\cA)/N}, \label{eq:scalingfixpoint}
\end{align}
with $G(\cA)/N = G(\cA_0)/H$. Here, $N$ is the counterpart of the charges that
disappear in the scaling limit, while $H$ corresponds to the confined charges
(i.e., those which appear only in the scaling limit).
In section~\ref{sec:coinciding} we show
\begin{equation}\label{eq:FN}
\big(\cF(\cA)^N\big)_0 = \cF(\cA)_0,
\end{equation}
which is again a case in which~\eqref{eq:equal} holds. Also, notice that here both nets involved satisfy (twisted) Haag duality.

In the remaining part of this work, we investigate the scaling limit of subsystems $\cA \subset \cB$ of the form $\cF^K \subset \cF^G$, where $\cF=\cF(\cA) =\cF(\cB)$ is a
graded-local field net acted upon by the compact groups $G \subset K$.
Of course, this situation includes the case of a field net and its gauge-invariant observable subnet recalled above,
but, for example, also the subnet generated by the local energy-momentum tensor fits in. In this framework,
we discuss the general relations between the groups appearing in equations~\eqref{eq:scalingfieldfixpoint}, \eqref{eq:scalingfixpoint} associated to $\cA$ and $\cB$. We then apply the results on classification of subsystems in~\cite{CC1, CC2} to gain some insights on the structure of the inclusion $\cA_0 \subset \cB_0$, and in particular on the relation between the canonical field nets in the scaling limit $\cF(\cA_0)$ and $\cF(\cB_0)$.

\medskip
The content of this paper is as follows.
In section 2
we show that for a subsystem $\cA \subset \cB$ with a conditional expectation $E : \cB \to \cA$ there is a one-to-one correspondence between the sets of scaling limit states of $\cA$ and those of $\cB$. This entails the somewhat curious fact that the sets of scaling limit states of
any two theories are in bijective correspondence~\cite{DM}. As another consequence, we show that the scaling limit of the $\bZ_2$-fixed point net of the free massive scalar field coincides with the $\bZ_2$-fixed point net of the free massless scalar field. Then
we readily adapt the argument in order to deal with more general free fields.
In section 3 we prove equation~\eqref{eq:FN}. Finally in section 4 we present a detailed discussion of the scaling limit of subsystems of the form $\cF^K \subset \cF^G$, illustrating the main results with several examples.

\section{Scaling limit for subsystems}\label{sec:subsystem}
We start by recalling some known facts to be used in the following, also to fix our terminology and notation.

\begin{Definition}
By a graded-local net with gauge symmetry we mean a quadruple $(\cF,\a,\b,\o)$, where:
\begin{proplist}{3}
\item $O \to \cF(O)$ is a net of unital C$^*$-algebras over double cones in Minkowski spacetime;
\item $\a$ is an automorphic action on $\cF$ of a geometrical symmetry group $\Gamma$ (the Poincar\'e group or its normal subgroup of translations) such that, for each double cone $O$, $\a_\gamma(\cF(O)) = \cF(\gamma\cdot O)$, $\gamma \in \Gamma$;
\item $\b$ is an automorphic action on $\cF$ of a compact group $G$ commuting with $\a$ and such that, for each double cone $O$, $\b_g(\cF(O)) = \cF(O)$, $g \in G$;
\item $\o$ is a pure state on $\cF$ which is $\a$- and $\b$-invariant;
\item there exists an element $k$ in the centre of $G$ with $k^2 = e$ such that, by defining
\begin{equation*}
F_\pm := \frac{1}{2}(F \pm\b_k(F)), \qquad F \in \cF,
\end{equation*}
for $F_i \in \cF(O_i)$, $i=1,2$, with $O_1$ spacelike from $O_2$, there holds
\begin{equation*}
F_{1,+}F_{2,\pm} = F_{2,\pm} F_{1,+}\,, \qquad F_{1,-}F_{2,-}= -F_{2,-}F_{1,-}\,.
\end{equation*}
\end{proplist}
\end{Definition}

We need also a spatial version of the above concepts.

\begin{Definition}
A graded-local net with gauge symmetry in the vacuum sector will be a graded local net with gauge symmetry such that:
\begin{proplist}{3}
\item
 for each $O$, $\cF(O)$ is a von Neumann algebra acting on the Hilbert space $\cH$;
\item there is a strongly continuous unitary representation $U$ of $\Gamma$ on $\cH$ such that $\a_\gamma = \ad \,U(\gamma)$, and such that the joint spectrum of the generators of the representation of the  translations subgroup $\bR^4 \ni x \to U(x)$ is contained in the closed forward light cone;
\item there is a strongly continuous unitary representation $V$ of $G$ on $\cH$ commuting with $U$ and such that $\b_g = \ad \,V(g)$;
\item $\o$ is the vector state induced by a $U$- and $V$-invariant unit vector $\O \in \cH$ which is cyclic for the quasi-local algebra $\cF = \overline{\bigcup_O \cF(O)}$ (closure in the norm topology).
\end{proplist}
\end{Definition}

If $G = \{e,k\}\cong \bZ_2$, we will simply speak of a graded-local net. In the particular case in which $G$ is trivial (and therefore $k = e$), we will use the traditional notation $\cA$ instead of $\cF$, and we will refer to the triple $(\cA, \a, \o^\cA)$ as a local net (in the vacuum sector if it applies). If $(\cF, \a,\b,\o^\cF)$ is a graded-local net with gauge symmetry, then one obtains a local net by defining
\begin{equation*}
\cA(O) := \cF(O)^G:=\{F \in \cF(O)\,:\, \b_g(F) = F, \,g \in G\}.
\end{equation*}
Moreover, an Haag dual net will be a local net in the vacuum sector such that $\cA(O) = \cA(O')'$, where as usual $\cA(O')$ is defined as the C$^*$-algebra generated by the $\cA(O_1)$ for all double cones $O_1 \subset O'$.

\medskip
We recall the construction of the scaling algebra of a graded-local net with gauge symmetry in the vacuum sector $\cF$~\cite{BV, DMV}: we consider the C$^*$-algebra of all bounded functions $\uF : \bR_+ \to \cF$, with norm $\norm{\uF} := \sup_{\l > 0}\norm{\uFl}$, endowed with the automorphic actions of $\Gamma$  and $G$ defined by
\begin{equation*}
\ua_\gamma(\uF)_\l := \a_{\gamma_\l}(\uFl), \qquad \ub_g(\uF)_\l := \b_g(\uFl),\qquad \gamma \in \Gamma,\, g \in G,\,\l > 0,
\end{equation*}
where $\gamma_\l = (\Lambda, \l x)$ if $\gamma = (\Lambda,x)$. Then $\uFF(O)$ is the C$^*$-subalgebra of the functions $\uF$ such that
\begin{enumerate}
\item $\uFl \in \cF(\l O)$ for all $\l > 0$;
\item $\displaystyle \lim_{\gamma \to (\id,0)} \norm{\ua_\gamma(\uF)-\uF}=0$;
\item $\displaystyle \lim_{g \to e}\norm{\ub_g(\uF)-\uF}=0$.
\end{enumerate}
In the particular case in which $\cF = \cA$ is a local net, the third condition above is of course void because of the triviality of $G$. We denote by $\uFF$ the quasi-local C$^*$-algebra defined by the net $O \to \uFF(O)$.

\begin{Remark}\label{rem:gfactor}
According to property 3.\ above, the scaling algebra $\uFF$ associated to $(\cF,\a,\b,\o)$ depends on
the action $\b$ of $G$. Since we do not require $\b$ to be faithful, it factors through an action $\tilde \b$ of $G/N$, where $N := \{ g \in G \,:\, \b_g(F) = F,\,\forall F \in \cF\}$, and one could consider the scaling algebra $\tilde \uFF$ associated to $(\cF,\a,\tilde \b,\o)$. However, thanks to the fact that the canonical projection $G \to G/N$ is open, it turns out straightforwardly that actually $\tilde \uFF = \uFF$.
\end{Remark}

Next we introduce scaling limits. We define states $\uo_\m$, $\m >0$, on $\uFF$ by $\uo_\m(\uF) := \o(\uF_\m)$, and we denote by $\textup{SL}(\o^\cF)$ the set of weak* limit points of $(\uo_\m)_{\m>0}$. We shall write $\textup{SL}(\o^\cF) = (\uooi)_{\iota \in \cI_\cF}$, where $\cI_\cF$ is an appropriate index set. Each $\uooi$ will be called a scaling limit state of $\o$,\footnote{A more general definition of scaling limit state has been given in~\cite{BDM}.} and we denote by $(\poi, \Hoi, \Ooi)$ the GNS triple induced by $\uooi$.  According to the results in~\cite[sec.\ 3]{DMV}, $(\uFF, \ua, \ub, \uooi)$ is a graded-local net with gauge symmetry, and by defining
\begin{equation}
\Foi(O) := \poi(\uFF(O))''
\end{equation}
one gets a graded-local net with gauge symmetry in the vacuum sector, called a scaling limit net of $\cF$. The notation $\ooi = \langle \Ooi,(\cdot)\Ooi\rangle$ will be systematically employed in the following.

\medskip
For a general analysis of the notion of subsystem see~\cite{Wi, LR, CDR}.

\begin{Definition}\label{def:inclusion}
Given two graded-local nets $(\cF, \a^\cF, \b^\cF, \o^\cF)$, $(\cB, \a^\cB, \b^\cB, \o^\cB)$, we say that they form an \emph{inclusion of graded-local nets}, and write for brevity $\cB \subset \cF$, if:
\begin{itemize}
\item[(i)] $\cB(O) \subset \cF(O)$ for each double cone $O$;
\item[(ii)] $\a_\gamma^\cF(B)= \a_\gamma^\cB(B)$, for all $B \in \cB$, $\gamma \in \Gamma$;
\item[(iii)] $\b_{k_\cF}^\cF(B)= \b_{k_\cB}^\cB(B)$, for all $B \in \cB$;  \label{it:inclusion:grading}
\item[(iv)]  $\o^\cF(B) = \o^\cB(B)$ for all $B \in \cB$.
\end{itemize}
\end{Definition}

Accordingly, when there is no danger of confusion, we will omit indices $\cF, \cB$ and write simply $\a$, $k$ and $\o$. In the above situation, if $(\cF, \a^\cF,\b^\cF,\o^\cF)$ is a graded-local net in the vacuum sector, it follows easily, by a Reeh-Schlieder type argument, that $\O$ is separating for $\cF(O)$ for each $O$, and therefore it is  clear that by restricting $\cB(O)$, $U_\cF(\gamma)$ and $V_\cF(k)$ to $\cH_\cB := \overline{\cB\O_\cF} \subset \cH_\cF$, one gets a graded-local net in the vacuum sector, which is isomorphic to $(\cB, \a^\cB, \b^\cB,\o^\cB)$ (see e.g.~\cite{CC1}, top of page 93), and therefore it will be identified with  $(\cB, \a^\cB, \b^\cB,\o^\cB)$ when no ambiguities arise.

\medskip

In the sequel,
we also assume the existence of a conditional expectation of nets $E:\cF\to\cB$, meaning that $E$ is a conditional expectation on the quasi-local algebra $\cF$ onto the quasi-local algebra $\cB$, which in restriction to every $\cF(O)$ is a conditional expectation onto $\cB(O)$, and such that $\a_\gamma E = E \a_\gamma$, $\b_k E = E \b_k$ and $\o\circ E = \o$. It follows from the last property that if $\cB$, $\cF$ are in the vacuum sector, $E$ restricts to a normal conditional expectation of $\cF(O)$ onto $\cB(O)$. Such a conditional expectation exists if, e.\ g., $\cF$ and $\cB$ satisfy twisted Haag duality on their respective vacuum spaces~\cite[sec. 3]{CDR} (see also~\cite{LR}). Our setup includes in particular the case where $\cF$ is a Doplicher-Roberts field net over a local net of observables $\cB$, so that $E$ is obtained by taking the average over the compact global gauge group~\cite{DR}.

\medskip
Now we wish to examine the possible relations between the scaling algebras $\uFF(O)$ and $\uBB(O)$ and the scaling limit states $\textup{SL}(\o^\cF) = (\uo^\cF_{0,\iota})_{\iota \in \cI_\cF}$ and $\textup{SL}(\o^\cB) = (\uo^\cB_{0,\iota})_{\iota \in \cI_\cB}$ associated to $\cF$ and $\cB$ respectively. It is clear that since $\cF$ and $\cB$ satisfy conditions (i)-(iii) of definition~\ref{def:inclusion}, the same is true for $\uFF$ and $\uBB$. It is then easy to see that the map $\uE$ defined on $\uFF$ by
\begin{equation}\label{eq:scalingexp}
\uE(\uF)_\l := E(\uFl), \qquad \uF \in \uFF, \, \l > 0
\end{equation}
is a conditional expectation of nets from $\uFF$ onto $\uBB$, commuting with $\ua$ and $\ub_k$. Moreover, if $E$ is faithful, then also $\uE$ is: if, for each $\l > 0$, $\uE(\uF^*\uF)_\l = E(\uFl^*\uFl) = 0$, then $\uFl =0$, i.e. $\uF=0$.

\begin{Proposition}
Let $\cB \subset \cF$ be an inclusion of graded-local nets and $E:\cF\to\cB$ a conditional expectation as before. Then $\textup{SL}(\o^\cF)=\textup{SL}(\o^\cB)\circ\uE$, and there is a bijective correspondence between $\cI_\cB$ and $\cI_\cF$ defined by mapping $\uo^\Boi \in \textup{SL}(\o^\cA)$, $\iota \in \cI_\cB$, to $\uo^\cB_{0,\iota}\circ\uE \in \textup{SL}(\o^\cF)$.
\end{Proposition}

\begin{proof}
Let $\uo^\cB_{0,\iota} \in \textup{SL}(\o^\cB)$, $\iota \in \cI_\cB$. Then, since $\uo^\cB_\m \circ \uE(\uF) = \o^\cB(E(\uF_\m)) = \o^\cF(\uF_\m) = \uo^\cF_\m(\uB)$, we have that $\uo^\cB_{0,\iota}\circ\uE \in \textup{SL}(\o^\cF)$. Also, if $\uo^\cB_{0,\iota}\circ\uE=\uo^\cB_{0,\k}\circ\uE$, then $\uo^\cB_{0,\iota}(\uB) = \uo^\cB_{0,\iota}\circ\uE(\uB)=\uo^\cB_{0,\k}\circ\uE(\uB)=\uo^\cB_{0,\k}(\uB)$ for all $\uB \in \uBB$, and the map defined in the statement is injective.

Conversely, let $\uo^\cF_{0,\iota} \in \textup{SL}(\o^\cF)$, $\iota \in \cI_\cF$. Then $\uo^\cF_{0,\iota}$ is a weak* limit point of $(\uo^\cF_\m)_{\m > 0}$, and therefore $\uo^\cB_{0,\iota}:=\uo^\cF_{0,\iota} \rest \uBB$ is a weak* limit point of $(\uo^\cF_\m\rest\uBB)_{\m > 0}$. But, for $\uB \in \uBB$, $\uo^\cF_\m(\uB) = \o^\cF(\uB_\m) = \o^\cB(\uB_\m) = \uo^\cB_\m(\uB)$, and then $\uo^\cB_{0,\iota}\in \textup{SL}(\o^\cB)$, so that $\uo^\cF_{0,\iota}= \uo^\cF_{0,\iota}\circ\uE = \uo^\cB_{0,\iota}\circ\uE$. This also shows that the above defined map is surjective, concluding the proof.
\end{proof}

As a consequence of the above proposition,  $\uE$ is a conditional expectation of the nets $(\uFF, \ua^\cF,\ub^\cF, \uo^\cB_{0,\iota}\circ\uE)$ and $(\uBB, \ua^\cB,\ub^\cB, \uo^\cB_{0,\iota})$. Also, denoting by $\poi^\cB$ and $\poi^\cF$ the scaling limit representations defined by $\uo^\cB_{0,\iota}$ and $\uo^\cF_{0,\iota}=\uo^\cB_{0,\iota}\circ\uE$ respectively, we see that $\poi^\cF$ is the representation induced from $\poi^\cB$ via $\uE$.

\begin{Remark} It follows form the previous result that if $\cB \subset \cF$, even without assuming the existence of a conditional expectation of $\cF$ onto $\cB$ the map $\uo^\cF_{0,\iota} \to \uo^\cF_{0,\iota}\rest\uBB$ induces a bijection between $\cI_\cF$ and $\cI_\cB$.\footnote{As it is clear from the proof, this does not really depend on the fact that $\cB \subset \cF$.} In order to see this, assume, for simplicity, that $\cB$ and $\cF$ are local nets, and consider, as in~\cite[prop.\ 3.5]{DM}, the tensor product theory $\cG := \cB\otimes\cF$, and the conditional expectations $E^\cB : \cG \to \cB\simeq\cB\otimes\bC\Id$, $E^\cF:\cG \to \cF\simeq\bC\Id\otimes\cF$ given respectively by $E^\cB(B\otimes F) = \o^\cF(F)B$, $E^\cF(B\otimes F) = \o^\cB(B)F$. According to the previous proposition, we have a bijection between $\cI_\cF$ and $\cI_{\cG}$ induced by $\uo^\cF_{0,\iota} \to \uo^\cF_{0,\iota}\circ\uE^\cF$, and a bijection between $\cI_{\cG}$ and $\cI_\cB$ induced by $\uo^{\cG}_{0,\iota} \to \uo^{\cG}_{0,\iota}\rest\uBB$, where, with a slight abuse, we identify $\uBB$ with the (isomorphic) subalgebra of $\underline{\fG}$ consisting of the functions $\l \to \uB_\l\otimes \Id$, $\uB \in \uBB$. It is then sufficient to show that $\uo^\cF_{0,\iota}\circ\uE^\cF\rest\uBB = \uo^\cF_{0,\iota}\rest \uBB$, but this follows at once from
\begin{equation*}
\uo^\cF_{0,\iota}(\uE^\cF(\uB)) = \lim_\k \o^\cF(E^\cF(\uB_{\l_\k}\otimes \Id)) = \lim_\k \o^\cB(\uB_{\l_\k})= \uo^\cF_{0,\iota}(\uB).
\end{equation*}
The case in which $\cB$ and $\cF$ are genuinely graded-local nets can be handled in a similar way up to the replacement of the tensor product and of the slice maps with their $\bZ_2$-graded versions.
\end{Remark}

In view of the above proposition, we indentify $\cI_\cB$ and $\cI_\cF$ and denote both simply by $\cI$, and for each $\iota \in \cI$ we denote by $\Boi$ and $\Foi$ the scaling limit nets obtained by the corresponding states in $\textup{SL}(\o^\cB)$ and $\textup{SL}(\o^\cF)$. Now one can show the existence of a conditional expectation on every scaling limit theory.

\begin{Proposition}\label{prop:Eoi}
Given an inclusion $\cB \subset \cF$ of graded-local nets in the vacuum sector, there is, for each $\iota \in \cI$, an inclusion $\Boi \subset \Foi$ of scaling limit nets. Furthermore, if a conditional expectation of nets $E:\cF\to\cB$ is given, there exists a conditional expectation of nets $\Eoi : \Foi \to \Boi$ uniquely defined by $\Eoi(\poi^\cF(\uF)) := \poi^\cF(\uE(\uF))$, $\uF \in \uFF$. Moreover, if $\eoi := [\Boi\Ooi]$, one has
\begin{equation}\label{eq:jonesproj}
\Eoi(F)\eoi = \eoi F\eoi, \qquad F \in \Foi.
\end{equation}
\end{Proposition}

\begin{proof}
Thanks to the Reeh-Schlieder property for $\Foi$, the net $\Boi$ is isomorphic to the net $O \to \poi^\cF(\uBB(O))'' \subset \Foi(O)$, which gives the inclusion $\Boi \subset \Foi$.

In order to show the existence of the conditional expectation $\Eoi$, we start by observing that, given $\uB_1, \uB_2 \in \uBB$, $\uF \in \uFF$ one has
\begin{equation*}\begin{split}
\big\langle\poi^\cF(\uB_1)\Ooi,&\poi^\cF(\uF)\poi^\cF(\uB_2)\Ooi\big\rangle = \uooi(\uB_1^*\uF\uB_2) \\
&= \uooi(\uE(\uB_1^*\uF\uB_2)) = \uooi(\uB_1^*\uE(\uF)\uB_2) \\
&= \big\langle\poi^\cB(\uB_1)\Ooi,\poi^\cB(\uE(\uF))\poi^\cB(\uB_2)\Ooi\big\rangle,
\end{split}\end{equation*}
which, taking into account the above mentioned isomorphism, shows that the map $\poi^\cF(\uF) \to \poi^\cF(\uE(\uF))$ is well-defined and $\sigma$-weakly continuous, and therefore extends uniquely to a $\sigma$-weakly continuous map $\Foi \to \Boi$ which is easily seen to be a conditional expectation of nets.

The properties $\aoi_\gamma\Eoi = \Eoi\aoi_\gamma$, $\boi_k \Eoi = \Eoi\boi_k$ and $\ooi \circ \Eoi = \ooi$ follow at once from the analogous properties for $\uE$.

In order to prove~\eqref{eq:jonesproj}, it is clear, by normality, that it is sufficient to prove  $\Eoi(\poi^\cF(\uF))\eoi = \eoi \poi^\cF(\uF)\eoi$, $\uF \in \uFF$. This is shown by choosing, for $\Phi \in \Hoi^\cF$, a sequence $B_n \in \poi^\cF(\uBB)$ such that $B_n\Ooi$ converges to $\eoi\Phi$, and then by evaluating
\begin{equation*}\begin{split}
\langle\Phi,\Eoi(\poi^\cF(\uF))\eoi\Phi\rangle &= \langle\eoi\Phi,\Eoi(\poi^\cF(\uF))\eoi\Phi\rangle\\
&= \lim_{n \to +\infty} \langle\Ooi,\Eoi(B_n^*\poi^\cF(\uF)B_n)\Ooi\rangle \\
&= \lim_{n \to +\infty} \langle\Ooi,B_n^*\poi^\cF(\uF)B_n\Ooi\rangle\\
&= \langle\eoi\Phi,\poi^\cF(\uF)\eoi\Phi\rangle = \langle\Phi,\eoi\poi^\cF(\uF)\eoi\Phi\rangle,
\end{split}\end{equation*}
which gives the statement.
\end{proof}

\begin{Remark}
The above discussion carries over to the case where $\cB \subset \cF$ is an inclusion of graded-local nets with gauge symmetry, meaning that condition (iii) in definition~\ref{def:inclusion} is replaced by the following
\begin{itemize}
\item[(iii')] $\b_g^\cB(B) = \b_{\phi(g)}^\cF(B)$, for all $B \in \cB$, $g \in G_\cB$,
\end{itemize}
where $\phi : G_\cB \to G_\cF$ is a continuous homomorphism such that $\phi(k_\cB) = k_\cF$, and there exists a conditional expectation of nets $E : \cF \to \cB$ such that $\b_g^\cB(E(F)) = E(\b_{\phi(g)}^\cF(F))$ for all $F \in \cF$, $g \in G_\cB$. Notice in fact that, if $\uF \in\uFF(O)$, then $\uE(\uF)$ defined by~\eqref{eq:scalingexp} is still an element of $\uBB(O)$, since
\begin{equation*}\begin{split}
\lim_{G_\cB\ni g \to e} \sup_{\l > 0}\| \b^\cB_g(E(\uF_\l)) - E(\uF_\l)\| &= \lim_{G_\cB\ni g \to e} \sup_{\l > 0}\| E(\b_{\phi(g)}^\cF(\uF_\l)) - E(\uF_\l)\| \\
&\leq  \lim_{G_\cB\ni g \to e} \sup_{\l > 0}\|\b_{\phi(g)}^\cF(\uF_\l) - \uF_\l\| = 0.
\end{split}\end{equation*}
For instance, this situation applies if $N \subset G_\cF$ is a closed normal subgroup, $\cB := \cF^N$ and $E$ is the average on $N$, in which case one can assume $G_\cB = G_\cF$ (and therefore $\phi = \textup{id}$), whose action factors through $G_\cF/N$.\footnote{We stress that, due to remark~\ref{rem:gfactor}, the scaling algebra of $\cB$ when it is thought with an action of $G_\cF$ coincides with the scaling algebra obtained when $\cB$ is thought with an action of $G_\cF/N$.} When $N = G_\cF$, $\cA := \cB$ is the local net of observables associated to $\cF$, and we recover the existence of a conditional expectation from $\Foi$ to the observable scaling limit net $\Aoi$ used in the proof of lemma 5.1 in~\cite{DMV}.
\end{Remark}

It is worth pointing out that our treatment is by no means limited to nets of von Neumann algebras. This should be already clear from the above discussion, and it is further exemplified by the next theorem
which, following closely the arguments expounded in \cite{BV2}, is presented in the setting of nets of $C^*$-algebras, and where, as an application of the above results, we compute the scaling limit of the even part of the free scalar field. A von Neumann  algebraic version will follow from the more general theorem~\ref{thm:fieldmultiplet} afterwards.

\medskip
We will use the description of the scalar field net in terms of time-zero fields and locally Fock representations as in~\cite{BV2}, as well as the main results of that paper. In particular, we consider the case where $\Gamma = \bR^d$, $d =3,4$. Also, we associate to the free scalar field of mass $m \geq 0$ the net of C$^*$-algebras $O\to\cFm(O)$ obtained by considering the elements of the canonical net of von Neumann algebras of the free scalar field of mass $m$ which are norm-continuous under translations. This net is conveniently isomorphically represented on the Fock space of the massless scalar field by local normality, see \cite{BV2} for details.

Since we deal with nets of C$^*$-algebras, in the following result it is understood that  the scaling limit of a net $\cF$ of C$^*$-algebras is defined as $\Foi(O) = \poi(\uFF(O))$, without weak closure. Furthermore, we denote by $\cF_r(O) := \cap_{O_1 \supset \bar{O}}\cF(O_1)$ the outer regularized net of $\cF$.

\begin{Theorem}\label{thm:free}
Let $(\cFm,\am,\om)$ be the net associated to the free neutral scalar field of mass $m \geq 0$ in $d = 3,4$ spacetime dimensions, as described above, and let $(\cAm,\am,\om)$ be the subnet of fixed points under the $\bZ_2$ action defined by the involutive automorphism
\begin{equation*}
\b(W(f)) = W(-f), \qquad f \in \cD(\bR^3).
\end{equation*}
Then each outer regularized scaling limit net $(\cAmoir,\amoi,\omoi)$ is isomorphic to the net $(\cAz,\az,\oz)$ of $\bZ_2$-fixed points of the net associated to the massless scalar field.
\end{Theorem}

\begin{proof}
In view of the above results, each net $(\cAmoir,\amoi,\omoi)$ is a subnet of some outer regularized scaling limit net $(\cFmoir,\amoi,\omoi)$. Let then $\phi: \cFmoir \to \cFz$ be the isomorphism onto the net of the massless scalar field, whose existence is proven in~\cite{BV2}. We will show that $\phi\rest\cAmoir$ is an isomorphism onto $\cAz$. We begin by showing that $\phi(\cAmoir) \subset \cAz$. To this end, let $\uA \in \uAm(O)$; then, since in the chosen representation $\b$ is weakly continuous, being implemented by the unitary operator $e^{i\pi N}$ (with $N$ the number operator of the massless scalar field), we have, for a suitable net $(\l_\k)_\k \cont \bR_+$,
\begin{equation*}\begin{split}
\b(\phi(\poi(\uA))) &= \text{w-}\lim_\k \b\s_{\l_\k^{-1}}(\uA_{\l_\k})\\
                    &= \text{w-}\lim_\k \s_{\l_\k^{-1}}\b(\uA_{\l_\k})\\
                    &= \text{w-}\lim_\k \s_{\l_\k^{-1}}(\uA_{\l_\k}) = \phi(\poi(\uA)),
\end{split}\end{equation*}
where, in the second equality, we have used the fact that $\b$ commutes with the dilations $\s_\l$ as defined in~\cite[eq.\ (2.6)]{BV2}, and in the third equality the fact that $\uA_\l$ is $\b$-invariant for each $\l$. Therefore $\phi(\poi(\uA)) \in \cAz(O) = \cFz(O)^{\bZ_2}$, and then
$$\phi(\cAmoir(O)) = \phi\bigg(\bigcap_{O_1 \supset \bar O}\cAmoi(O_1)\bigg)\subset \cAz(O)$$
thanks to the outer regularity of $\cAz(O)$.

Conversely, $\phi$ being an isomorphism, any element $A \in \cAz(O)$ is of the form $A = \phi(\poi(\uF_1))$ for any $O_1 \supset \bar O$ and for some $\uF_1 \in \uFm(O_1)$. For such an element, define $\uA_1 := \uE(\uF_1) \in \uAm(O_1)$, where $E = \frac{1}{2}(\textup{id}+\b)$ is a conditional expectation of $\cFm$ onto $\cAm$, which is obviously weakly continuous and commuting with the dilations. Then, arguing as above, we have
\begin{equation*}\begin{split}
\phi(\poi(\uA_1)) &= \text{w-}\lim_\k \s_{\l_\k^{-1}}(E(\uF_{1,\l_\k})) \\
                &= E(\phi(\poi(\uF_1))) = \phi(\poi(\uF_1)) = A,
\end{split}\end{equation*}
so that $\phi(\cAmoir) = \cAz$.
\end{proof}

\begin{Remark}
The net $\cA^{(m_1)}\otimes \cA^{(m_2)}$ with $m_1 \neq m_2$ has no
nontrivial
subsystems thanks to the results in~\cite{CC1}, but, according to the above theorem and the results in~\cite{DM}, any of its scaling limits is  isomorphic to $\cA^{(0)}\otimes \cA^{(0)}$, which has, for instance, the subsystem obtained by taking the fixpoint net with respect to the natural action of $SO(2)$. Therefore, this simple example shows that subsystems can appear in the scaling limit which are not related to subsystems already existing at finite scales.
\end{Remark}

As anticipated in~\cite[thm.\ 4.6]{DM} the argument in theorem~\ref{thm:free} carries over to the case of multiplets of free fields acted upon by a compact Lie group $G$. More precisely, consider a finite symmetric and generating set $\D$ of irreducible representations of $G$ and a mass function $\bmu : \D \to [0,+\infty)$ such that $\bmu(v) = \bmu(\bar{v})$. Let then $\cFbm$ denote the graded-local net with gauge symmetry in the vacuum sector generated by a $v$-multiplet of free scalar fields of mass $\bmu(v)$ for each $v \in \D$, and  $\cAbm \subset \cFbm$ the fixed point net of $\cFbm$ under the natural action of $G$. Again, this net is isomorphically represented on the Fock space $\cHbz$ corresponding to $\bmu(v) =0$ for each $v \in \D$ (see~\cite{DM} for details). Furthermore, denote by $\cAbmoi \subset \cFbmoi$ the corresponding inclusion of scaling limit nets. As shown in~\cite[thm 4.3]{DM} there is a spatial net isomorphism $\theta : \cFbmoi \to \cFbz$ such that for each $\uF \in \uFF^{(\bmu)}$,
\begin{equation*}
\theta(\pbmoi(\uF)) = \text{w-}\lim_{\k}\d^{-1}_{\l_\k}(\uF_{\l_\k}),
\end{equation*}
for a suitable net $(\l_\k)_\k \subset \bR_+$, and where $\d_\l$ is the adjoint action of the dilation group on $\cHbz$.

\begin{Theorem}\label{thm:fieldmultiplet}
There is a net isomorphism between $\cAbmoi$ and $\cAbz$, obtained from $\theta$ by restriction.
\end{Theorem}

\begin{proof}
Since the action of $G$ on $\cFbm$ is $\bmu$-independent~\cite{DM}, the same is true for the conditional expectation $E : \cFbm \to \cAbm$ obtained by averaging with respect to $G$. Therefore, if $\Ebmoi : \cFbmoi \to \cAbmoi$ is the conditional expectation given by proposition~\ref{prop:Eoi}, in order to generalize the above argument it is sufficient to show that $\theta \circ \Ebmoi = E \circ\theta$. By normality of $\theta$ and $\Ebmoi$, it is sufficient to check this equation on elements $\pbmoi(\uF)$ with $\uF \in \uFF^{(\bmu)}(O)$ for some $O$. This follows at once from the computation
\begin{equation*}\begin{split}
\theta\circ\Ebmoi(\pbmoi(\uF)) &= \text{w-}\lim_{\k}\d^{-1}_{\l_\k}(E(\uF_{\l_\k})) \\
                                                      &=\text{w-}\lim_{\k}E(\d^{-1}_{\l_\k}(\uF_{\l_\k})) = E\circ\theta(\pbmoi(\uF)),
\end{split}\end{equation*}
where in the last equality the norm boundedness of $\d^{-1}_{\l_\k}(\uF_{\l_\k})$ and the normality of $E$ were used.
\end{proof}

\medskip
The essential point in the above proofs is the existence of conditional expectations $E_{0,\iota}^{(\bmu)}: \cFbmoi \to \cAbmoi$ and $E : \cFbz \to \cAbz$ intertwining the action of the isomorphism $\theta : \cFbmoi \to \cFbz$. In fact if we consider the general situation of an inclusion $\cB \subset \cF$ with a conditional expectation of nets $E:\cF\to\cB$ as discussed above, and we assume that $\cF$ has a unique (quantum) scaling limit, with isomorphisms $\phi_{\iota,\iota'}: \Foi \to \cF_{0,\iota'}$, and that the conditional expectations $\Eoi : \Foi \to \Boi$ introduced in proposition~\ref{prop:Eoi} satisfy
\begin{equation}\label{eq:conditionalinter}
\phi_{\iota,\iota'}\circ\Eoi = E_{0,\iota'}\circ\phi_{\iota,\iota'},
\end{equation}
a similar argument shows that $\phi_{\iota,\iota'}(\Boi) = \cB_{0,\iota'}$, so that $\cB$ has a unique scaling limit too.

This happens in particular if $\cF$ has a \emph{convergent scaling limit} as introduced in~\cite{BDM}:  we say that a net $\cF$ has a convergent scaling limit if there exists an inclusion of graded-local nets $\hat \uFF \subset \uFF$ such that for each $\uF \in \hat\uFF$ there exists $\lim_{\l \to 0}\o(\uF_\l)$ and such that, for each scaling limit state $\uooi$, in the corresponding scaling limit representation  one has $\poi(\hat\uFF(O))'' = \Foi(O)$.  It is easily seen that if a theory has a convergent scaling limit then the scaling limit is unique.

\begin{Proposition}
Let $\cB \subset \cF$ be an inclusion of graded-local nets in the vacuum sector, such that $\cF$ has convergent scaling limit, and let $E : \cF \to \cB$ be a conditional expectation of nets. Then equation~\eqref{eq:conditionalinter} holds. Furthermore $\cB$ has convergent scaling limit too.
\end{Proposition}

\begin{proof}
It is straightforward to show that the unitary $V_{\iota,\iota'} : \Hoi \to \cH_{0,\iota'}$ defined by
\begin{equation*}
V_{\iota,\iota'}\poi(\hat\uF)\Ooi = \pi_{0,\iota'}(\hat \uF)\O_{0,\iota'}, \qquad \hat \uF \in \hat \uFF,
\end{equation*}
implements a net isomorphism  $\phi_{\iota,\iota'}: \Foi \to \cF_{0,\iota'}$, and there holds, for each $\hat \uF \in \hat \uFF(O)$,
\begin{equation*}
\phi_{\iota,\iota'}\circ\Eoi(\poi(\hat \uF)) = \pi_{0,\iota'}(\uE(\hat \uF)) = E_{0,\iota'}\circ\phi_{\iota,\iota'}(\poi(\hat\uF)),
\end{equation*}
so that equation~\eqref{eq:conditionalinter} follows thanks to $\poi(\hat\uFF(O))'' = \Foi(O)$.

Furthermore, using the normality of $\Eoi$, it is direct to verify that $\cB$ has a convergent scaling limit by setting $\hat \uBB(O) := \uE(\hat\uFF(O))$.
\end{proof}

\section{Inclusions with coinciding scaling limits}\label{sec:coinciding}
In the previous section we have discussed situations in which an inclusion of nets gives rise to a proper inclusion of nets in the scaling limit. For completeness, in the present section we provide a general construction of an inclusion of nets such that the corresponding inclusion of scaling limit nets is trivial.

\medskip
Let $(\cF,\a,\b,\o)$ be a graded-local net with gauge symmetry in the vacuum sector. Then the quintuple $(\cF, U,V,\O,k)$ is a QFTGA according to \cite{DMV}, where the representations $\Uoi$ and $\Voi'$ of the translations and of $G$ for the scaling limit $\Foi$ are introduced. As $\Voi'$ is not necessarily faithful, we define the closed normal subgroup
\begin{equation}\label{eq:N}
N := \{g \in G\,:\,\Voi'(g) = I\}.
\end{equation}
The scaling limit net $\Foi$ is then obviously covariant with respect to the natural representation $\Voi$ of the factor group $\Goi := G/N$.

\begin{Proposition}With the above notation, let $\cB$ be the subsystem of fixed points of $\cF$ under $N$, with its natural action of $G$. Then for the associated scaling limit net there holds $\Boi = \Foi$.
\end{Proposition}

\begin{proof}
From $\uBB(O) \subseteq \uFF(O)$, $\Boi(O) \subseteq \Foi(O)$ readily follows. In order to prove the reverse inclusion, take $\uF \in \uFF(O)$ and define
\begin{equation*}
\uB := \int_N dn\, \ub_n(\uF),
\end{equation*}
where the integral is performed with respect to the normalized Haar measure on $N$ and is understood in Bochner sense, cfr.~\cite{Y}. This is well defined, since, by definition of $\uFF(O)$, $n \in N \to \ub_n(\uF) \in \uFF(O)$ is a continous function on a compact space, and then its range, being metrizable, is separable. We obtain then that $\uB \in \uFF(O)$. Furthermore, as the function $\uF \in \uFF(O) \to \uF_\l \in \cF(\l O)$ is norm continuous, and a Bochner integral is a norm limit of Lebesgue sums, we get
\begin{equation*}
\uB_\l = \bigg(\int_N  dn\, \ub_n(\uF)\bigg)_\l = \int_N dn\, \ub_n(\uF)_\l = \int_N dn\, \b_n(\uF_\l)
\end{equation*}
so that, for $m \in N$,
\begin{equation*}
\b_m(\uB_\l) = \int_Ndn\, \b_{mn}(\uF_\l) = \uB_\l,
\end{equation*}
having used the invariance of the measure $dn$. This shows then that $\uB \in \uBB(O)$. Then, using again the norm continuity of $n \to \ub_n(\uF)$, that of $\poi$, and the definition of $\Voi'$, we get
\begin{equation*}
\poi(\uB) = \int_N dn\, \poi(\ub_n(\uF)) = \int_N dn\, \Voi'(n) \poi(\uF)\Voi'(n)^* = \poi(\uF),
\end{equation*}
where the last equality follows from the definition of $N$. Thus we get $\poi(\uF) \in \Boi(O)$ and the statement of the proposition.
\end{proof}

\begin{Remarks}
\begin{remlist}
\item At first sight, one might think that the above result is a trivial consequence of lemma 5.1.(i) of~\cite{DMV}, but some subtleties in the definition of the relevant scaling algebras prevent the application of the cited result. This is because the definition of the scaling limit net $\Foi$ really depends on \emph{both} the underlying net $\cF$ \emph{and} the group $G$ acting on it, through the requirement of norm continuity of functions $g \in G \to \ub_g(\uF)$, $\uF \in \uFF$. So, if one was willing to apply lemma 5.1.(i) of~\cite{DMV} with $\cA = \cB$, he should define a \emph{new} scaling net $\tilde{\uFF}$ associated to the datas $(\cF,N)$, i.e. in the same way as $\uFF$ but requiring now only continuity of $n \in N \to \ub_n(\uF)$. In general, this would result in a much bigger net than $\uFF$. Then, application of the cited result would lead to $\Boi(O) = \tilde{\cF}_{0,\iota}(O)^{N/\tilde{N}}$, where now $\tilde{N}$ is a normal subgroup of $N$, defined in the obvious way. Also, the scaling limit net would now be acting on a new Hilbert space, in general much bigger than our $\Hoi$. There are however cases in which the scaling net $\uFF$ does not really depend on the group $G$, and then the result of~\cite{DMV} can be applied straightforwardly. For instance, this is the case if $G$ is a finite group, so that the continuity requirement is empty, which entails $\tilde{\uFF} = \uFF$, $\tilde{N} = N$ and finally $\Boi(O) = \cF_{0,\iota}(O)^{N/N} = \cF_{0,\iota}(O)$.
\item The group $N$ is really non-trivial, in general: if $\phi_i$, $i=1,\dots,n$, are charged generalized free scalar fields with mass measure $d\rho(m^2) = c\, dm^2$, on which a compact gauge group $G \subseteq U(n)$ acts, and if $\cF(O)$ is generated by the fields $\Box^{n(O)}\phi_i(f)$ with $\supp f \subset O$, where $n(O) \to +\infty$ as the radius of $O$ shrinks to $0$, then the scaling limit net $\cF_0$ is trivial~\cite{DM}, and therefore $\Voi'(g)=I$ for all $g \in G$, and $N=G$. More generally, in~\cite{DM} examples are constructed where $N$ is any closed subgroup of an arbitrary compact Lie group $G$.
\item  For any net $O \to \cC(O)$ such that $\cB(O) \subseteq \cC(O) \subseteq \cF(O)$, we define the associated ``interpolated'' scaling algebras as
\begin{equation*}
\underline{\fC}(O) := \{ \uF \in \uFF(O) \,:\, \uF_\l \in \cC(\l O)\},
\end{equation*}
and the corresponding scaling limit net as
\begin{equation*}
\cC_{0,\iota}(O) := \poi(\underline{\fC}(O))''.
\end{equation*}
Then it follows at once from the above proposition that $\cC_{0,\iota}(O) = \Foi(O)$.
\end{remlist}
\end{Remarks}

\section{Scaling of subsystems and Doplicher-Roberts reconstruction}\label{sec:DR}
Inside a net of local observables, there are operators with a specific physical interpretation like the energy momentum tensor, or Noether currents associated to (central) gauge symmetries, and the relations between the given net and the subsystem generated by such operators have been thoroughly investigated in~\cite{CDR, CC1,CC2} from the point of view of Doplicher-Roberts (DR) theory. In the present context, it is therefore natural to analyse the scaling limit of such subsystems and characterize them as subsystems of the scaling limit.

\subsection{General properties}\label{subsec:general}
As a first step in this direction, in this section we deal with the following abstract situation: we consider an inclusion $\cA \subset \cB$ of Haag dual and Poincar\'e covariant nets in the vacuum sector as defined in section~\ref{sec:subsystem}. We also require that the vacuum Hilbert space $\cH_\cB$ is separable. Thanks to the results in the appendix of~\cite{Rob} (see also the remark in sec.\ 4 of~\cite{DMV}), the main results in~\cite{DR} can be applied to $\cA$ and $\cB$, and we further assume that for the corresponding DR canonical  covariant field nets one has $\cF(\cA) = \cF(\cB)$.  Therefore for the canonical DR gauge groups one has that $G(\cB)$ is a closed subgroup of $G(\cA)$.

 Let us fix a scaling limit state $\uo^\cB_{0,\iota}$ of $\cB$. According to the results in~\cite{DMV}, there exists a scaling limit $\cF(\cB)_{0,\iota}$ of $\cF(\cB)$ and a quotient $G(\cB)_{0,\iota}$ of $G(\cB)$ by a normal closed subgroup $N(\cB)_{0,\iota}$ defined in analogy to~\eqref{eq:N}, such that $\Boi = \cF(\cB)_{0,\iota}^{G(\cB)_{0,\iota}}$. Furthermore if $\Boi$ satisfies Haag duality and if its vacuum Hilbert space is separable, denoting by $\cF(\Boi)$ the canonical DR field net of $\Boi$, one has that $\cF(\cB)_{0,\iota}$ is a fixed point net of $\cF(\Boi)$ with respect to a certain normal closed subgroup $H(\Boi)$ of the canonical DR group $G(\Boi)$. Thanks to what was shown in section~\ref{sec:subsystem}, using the corresponding scaling limit state $\uo^\cA_{0,\iota} := \uo^\cB_{0,\iota}\rest\uAA$ of $\cA$ we get similar relations for the nets $\Aoi$, $\cF(\cA)_{0,\iota}$ and $\cF(\Aoi)$. Summarizing, we get the following result.

\begin{Proposition}
With the above notations, the following diagram of inclusions of nets holds:
\begin{equation}\label{eq:diagram}
\begin{array}{ccccc}
\Boi &\subset &\cF(\cB)_{0,\iota} &\subset &\cF(\Boi) \\
\cup          &        &\cup               &        &\cup \\
\Aoi &\subset &\cF(\cA)_{0,\iota} &\subset &\cF(\Aoi)
\end{array}
\end{equation}
\end{Proposition}

\begin{proof}
As noted above, the horizontal lines follow from the results in~\cite{DMV}, while the first column is a consequence of the discussion in section~\ref{sec:subsystem}. The second column is immediate from the definition of the scaling limit net and the fact that $G(\cB) \subset G(\cA)$, and finally the third column follows from the first and~\cite{CDR}.
\end{proof}

 Notice that, even if $\cF(\cA) = \cF(\cB)$, the results of~\cite{DMV} do not allow to conclude that $\cF(\cA)_{0,\iota}= \cF(\cB)_{0,\iota}$ because of the fact that the construction of $\cF(\cB)_{0,\iota}$ depends on $G(\cB)$ (and similarly for $\cF(\cA)_{0,\iota}$), see~\cite[def.\ 2.2]{DMV}.

\medskip
For completeness we also analyse the relations between the different gauge groups that arise in
diagram~\eqref{eq:diagram}. According to~\cite[sec.\ 2, 5]{DMV} and to the previous discussion, we have groups $G(\cA)$, $N(\cA)_{0,\iota}$,
$G(\cA)_{0,\iota}$, $G(\Aoi)$ and $H(\Aoi)$ such that
$G(\cA)/N(\cA)_{0,\iota} = G(\cA)_{0,\iota} = G(\Aoi)/H(\Aoi)$, and similarly for $\cB$.

\begin{Theorem}
Under the standing assumptions, we have that $N(\cB)_{0,\iota}$ is a subgroup of $N(\cA)_{0,\iota}$ and that there exists a morphism $\phi : G(\Boi) \to G(\Aoi)$ such that $\phi(H(\Boi))\subset H(\Aoi)$, and such that the quotient morphism $\tilde{\phi}$ on $G(\Boi)/H(\Boi) = G(\cB)/N(\cB)_{0,\iota}$ is given by $\tilde{\phi}(gN(\cB)_{0,\iota})=gN(\cA)_{0,\iota}$. Moreover, if $\cF(\Boi)=\cF(\Aoi)$, then $\phi$ is injective, and if in addition $\cF(\cB)_{0,\iota} = \cF(\cA)_{0,\iota}$, then $N(\cB)_{0,\iota} = N(\cA)_{0,\iota}\cap G(\cB)$, $H(\Boi) = H(\Aoi)$, and $\tilde{\phi}$ is injective too.
\end{Theorem}

\begin{proof}
As already remarked, we have $\uFF(\cA)\subset\uFF(\cB)$ and that $G(\cB)$ is a subgroup of $G(\cA)$. If $\uo_{0,\iota} = \lim_\k \uo_{\l_\k}$, it is easy to see that
\begin{equation*}
N(\cA)_{0,\iota} = \Big\{ g \in G(\cA)\,:\,\lim_\kappa \| (\beta_g(\uF_{\l_\k})-\uF_{\l_\k})\O\| = 0,\, \forall \uF \in \bigcup_{O}\uFF(\cA)(O)\Big\},
\end{equation*}
which immediately implies the inclusion $N(\cB)_{0,\iota}\subset N(\cA)_{0,\iota}$. The existence of the morphism $\phi$, as well as its injectivity in the case $\cF(\Boi)=\cF(\Aoi)$, are direct consequences of the application of~\cite[thm.\ 2.3]{CDR} to the commuting square of inclusions provided by $\Aoi$, $\Boi$, $\cF(\Aoi)$ and $\cF(\Boi)$. Since $\phi$ is given by the restriction to $\cF(\Aoi)$ of automorphisms of $\cF(\Boi)$, and $H(\Boi)$ is the subgroup of $G(\Boi)$ which leaves $\cF(\cB)_{0,\iota}$ pointwise invariant, and similarly for $H(\Aoi)$, it is clear that $\phi(H(\Boi))\subset H(\Aoi)$. This also implies that $\tilde{\phi}$ is the restriction to $\cF(\cA)_{0,\iota}$ of automorphisms of $\cF(\cB)_{0,\iota}$, and therefore it coincides with the map $gN(\cB)_{0,\iota} \to gN(\cA)_{0,\iota}$.

We assume now that $\cF(\cB)_{0,\iota} = \cF(\cA)_{0,\iota}$ and $\cF(\Boi) = \cF(\Aoi)$. It follows immediately that
$H(\Boi) = H(\Aoi)$. We show that $N(\cA)_{0,\iota} \cap G(\cB) \subset N(\cB)_{0,\iota}$, the reverse inclusion being trivial. Let $g \in N(\cA)_{0,\iota} \cap G(\cB)$, i.e.\ $g \in G(\cB)$ and
\begin{equation*}
\lim_\kappa \| (\beta_g(\uF_{\l_\k})-\uF_{\l_\k})\O\| = 0,\quad \forall \uF \in \bigcup_{O}\uFF(\cA)(O).
\end{equation*}
Then, if $\uF' \in \uFF(\cB)(O)$ we can find a norm-bounded sequence $\uF_n \in \uFF(\cA)(O)$ such that $\poi(\uF_n)$ converges
strongly to $\poi(\uF')$ as $n \to +\infty$. We have then
\begin{equation*}
\lim_\kappa \| (\beta_g(\uF'_{\l_\k})-\uF'_{\l_\k})\O\| \leq \| \poi(\ub_g(\uF'-\uF_n))\O_0\|+ \| \poi(\uF'-\uF_n)\O_0\|, \end{equation*}
which, together with the fact that $\ub_g$ is unitarily implemented in $\poi$, readily gives $g \in N(\cB)_{0,\iota}$.The injectivity of $\tilde{\phi}$ then clearly follows from $N(\cA)_{0,\iota} \cap G(\cB) = N(\cB)_{0,\iota}$.
\end{proof}

\subsection{Field nets with trivial superselection structure in the scaling limit}
Until now we have employed the minimal set of assumptions on the scaling limit nets which allow us to make sense of the elements in diagram~\eqref{eq:diagram}. In order to proceed further in the discussion of its properties, it is useful at this point to apply the general machinery that has recently become available in the theory of subsystems, which requires some rather natural additional assumptions on the scaling limit nets, see definition~\ref{def:trivialsuper}.
 Partial results on the problem of deriving such assumptions from suitable hypotheses on the underlying nets at scale $\l = 1$ are discussed below in this section. We hope to give a more thorough analysis of these issues somewhere else.

\medskip
Below, we present some examples which corroborate the natural conjecture that, at least in typical cases, we have an equality in the last column of diagram~\eqref{eq:diagram}.
In subsection~\ref{subsec:gmax}, we outline a strategy for proving that $\cF(\Boi) = \cF(\Aoi)$.
Thanks to theorem~\ref{thm:Q}, the main point will be to show that $\Aoi = \cF(\cB)_{0,\iota}^{G(\cA)}$.

\begin{example}\label{ex:free}
Let $\cB$ be the net generated by a $G$-multiplet of massive free scalar fields. Then $\cF(\cB) = \cB$ and $G(\cB)$ is trivial~\cite{Dri}. Let also $\cA = \cB^G$, so that $\cF(\cA)=\cB=\cF(\cB)$ and $G(\cA) = G$. From the arguments in~\cite{DM} it is possible to prove that, for each scaling limit state of $\cB$, $\cF(\cA)_{0,\iota} = \cF(\cB)_{0,\iota}$, and therefore diagram~\eqref{eq:diagram} trivially reduces to
$$
\begin{array}{ccccc}
\Boi &=            &\cF(\cB)_{0,\iota} &= &\cF(\Boi) \\
\cup &             &\shortparallel        &   &\shortparallel \\
\Aoi &\subset &\cF(\cA)_{0,\iota} &= &\cF(\Aoi)
\end{array}
$$
The equality $\cF(\cA)_{0,\iota} = \cF(\cB)_{0,\iota}$ is obtained in the following way: one first observes that $\cF(\cB)_{0,\iota} \simeq \cB^{(0)}$, the net generated by a corresponding $G$-multiplet of massless free scalar fields transforming under the same representation~\cite{BV2,DM}. We recall that in the non-standard free field representation used in~\cite{DM} (see also sec.~\ref{sec:subsystem}), for each double cone $O$ based on the time zero plane one has $\cB^{(0)}(O) = \cB(O)$. Then, for each such double cone $O$ and each DR $G$-multiplet $\psi_j \in \cB^{(0)}(O)$, we define, for a continuous compactly supported function $h$ on $\bR^4$,
\begin{equation*}
(\ua_h \psi_j)_\l := \int_{\bR^4}dx\,h(x) \a_{\l x}\d_\l(\psi_j),
\end{equation*}
so that $\ua_h\psi_j \in \uFF(\cA)(O_1)$ for suitable $O_1 \supset O$. One then shows, using the same arguments as in the proof of~\cite[thm.\ 3.1]{BV2}, that $\pi_{0,\iota}(\ua_h\psi_j) \in \cF(\cA)_{0,\iota}(O_1)$ converges strongly to $\psi_j$ as $h \to \d$, and therefore, by outer regularity, $\psi_j\in \cF(\cA)_{0,\iota}(O)$, which entails $\cF(\cB)_{0,\iota} \subset \cF(\cA)_{0,\iota}$. The converse inclusion being trivial, the conclusion follows.
\end{example}

\begin{example}
The equality $\cF(\Boi) = \cF(\Aoi)$, which holds in the above example, can be deduced
under suitable assumptions
from the fact that $\cF(\cB)_{0,\iota} = \cF(\cA)_{0,\iota}$, as shown e.g.\ by theorem~\ref{thm:Q} and remark~\ref{rem:Q}.
The latter condition is trivially satisfied if for instance $G(\cB)$ is open in $G(\cA)$,
or if $[G(\cA):G(\cB)]$ is finite.
A discussion of more general conditions under which this is true seems to be of independent interest but for the time being it will be postponed.
\end{example}

\begin{example}
Suppose that $\cF$ is a dilatation covariant graded-local net satisfying the Haag-Swieca compactness condition. Since $\cF$ is considered to have a trivial gauge group, it is net-isomorphic to any of its scaling limit $\cF_{0,\iota}$ through
\begin{equation}\label{eq:dilaiso}
\phi(\pi_{0,\iota}(\uF))= \textup{s-}\lim_\k \d_{\l_\k}^{-1}(\uF_{\l_\k}),
\end{equation}
where $(\d_\l)_{\l>0}$ are dilatations on $\cF$, see~\cite[prop.\ 5.1]{BV} (in this reference only observable nets are considered, but the generalization to nets having normal commutations relations is not difficult). Assume now that $\cF = \cF(\cA)$ is obtained as the DR field net of an observable net $\cA$, with gauge group $G=G(\cA)$. Of course the scaling limit $\cF(\cA)_{0,\iota}$ of $\cF(\cA)$ satisfies $\cF(\cA)_{0,\iota} \subset \cF_{0,\iota}$. We show that the converse inclusion also holds. It suffices to show that $\poi(\uFF(O))\subset \poi\big(\uFF(\cA)(O)\big)''$, where $\poi$ is, as before, the scaling limit representation of $\uFF$. Consider then $\uF \in \uFF(O)$ and $F = \phi(\pi_{0,\iota}(\uF)) \in \cF(O)$. Defining $\tilde{F} \equiv \b_\psi(F) := \int_G dg\, \psi(g)\b_g(F)$, where $\psi\in C(G)$, and $\tilde{\uF}_\l = \d_\l(\tilde F)$, we have $\tilde{\uF} \in \uFF(\cA)(O)$ and $\phi(\pi_{0,\iota}(\tilde{\uF})) = \tilde{F}$ converges strongly to $F$ as $\psi \to \d$. Therefore, $\phi$ being spatial, we obtain the desired inclusion and $\cF(\cA)_{0,\iota}=\Foi$.
Finally, if we also assume that $\cF= \cF(\cA) = \cF(\cB)$, with $\cA \subset \cB$, it is possible to show that $\cF(\cA)_{0,\iota} = \cF(\cB)_{0,\iota} = \cF(\Boi) = \cF(\Aoi) \simeq \cF$. To this end, we notice that the isomorphism $\phi_\cA$ between $\cF(\cA)_{0,\iota}$ and $\cF$, defined as in~\eqref{eq:dilaiso}, is just the restriction to $\cF(\cA)_{0,\iota}$ of the analogous isomorphism $\phi_\cB$ between $\cF(\cB)_{0,\iota}$ and $\cF$. Thus $\cF(\cA)_{0,\iota} = \cF(\cB)_{0,\iota} \simeq \cF$ and the conclusion follows from the isomorphisms $\cA \simeq \Aoi$, $\cB \simeq \Boi$. Therefore diagram~\eqref{eq:diagram} becomes in this case
$$
\begin{array}{ccccc}
\Boi &\subset       &\cF(\cB)_{0,\iota} &= &\cF(\Boi) \\
\cup &                   &\shortparallel        &   &\shortparallel \\
\Aoi &\subset       &\cF(\cA)_{0,\iota} &= &\cF(\Aoi)
\end{array}
$$
\end{example}

In the remaining part of this section we give a closer look at the situation in which for the nets of von Neumann algebras in the scaling limit there holds $\cF(\cA)_{0,\iota}=\cF(\cB)_{0,\iota}$.
Actually, we discuss the seemingly more general case in which  $\Aoi$ is the fixpoint net of $\cF(\cB)_{0,\iota}$ under a compact group action.

\begin{Definition}\label{def:trivialsuper}
We say that a graded-local net with gauge symmetry in the vacuum sector $(\cF,\a,\b,\O)$ has \emph{trivial superselection structure} if
\begin{proplist}{3}
\item $\O$ is cyclic and separating for $\cF(O)$ for each $O$ (Reeh-Schlieder property);
\item with $Z := (I+i k)/(1+i)$, there holds $\cF(O')'' = Z\cF(O)'Z^*$ for each $O$ (twisted Haag duality);
\item if $(\D,J)$ are the modular objects associated to $(\cF(W_R)'',\O)$, $W_R$ being the right wedge, there holds $\D^{it}=U(\L_{W_R}(t))$, $JU(\L,a)J = U(j\L j,j a)$ and, for each $O$, $J\cF(O)J = \cF(jO)$, where $\L_{W_R}$ is the one parameter group of boost leaving $W_R$ invariant and $j$ is the reflection with respect to the edge of $W_R$ (geometric modular action);
\item for each pair of double cones $O_1, O_2$ with $\bar{O}_1 \subset O_2$ there exists a type I factor $\cN_{O_1,O_2}$ such that $\cF(O_1) \subset \cN_{O_1,O_2}\subset \cF(O_2)$ (split property);
\item there exists at most one fermionic irreducible DHR representation with finite statistics $\pi$ of $\cF^{\bZ_2}$ inequivalent to the vacuum, and such that every DHR representation of $\cF^{\bZ_2}$ is the direct sum of copies of the vacuum representation and of $\pi$ (triviality of the superselection structure).
\end{proplist}
\end{Definition}

For a discussion of the above properties, in particular the last one, we refer the reader to~\cite{CC1,CC2}.

\begin{Theorem}\label{thm:Q}
Let $\cA \subset \cB$ be an inclusion of local nets satisfying the standing assumptions, and suppose furthermore that $\cF(\Boi)$ satisfies the properties (i)-(v) in definition \ref{def:trivialsuper}, and that $\Aoi = \cF(\cB)_{0,\iota}^Q$ for some compact group $Q$ of internal symmetries of the net $\cF(\cB)_{0,\iota}$. Then $\cF(\Aoi) = \cF(\Boi)$, and therefore $\Aoi$ is the fixpoint net of $\cF(\Boi)$ under a compact group of internal symmetries.
\end{Theorem}

\begin{proof}
As a first step, we show that $\cF(\cB)_{0,\iota} \subset \cF(\Aoi)$. To this end, we observe that, since $\Aoi = \cF(\cB)_{0,\iota}^Q$, $\cF(\cB)_{0,\iota}$ inside $\cF(\Boi)$ is generated (by $\Aoi$ and) by the Hilbert spaces of isometries in $\cF(\Boi)$ implementing the cohomological extensions to $\Boi$ of the covariant DHR sectors of $\Aoi$ corresponding to the irreducible representations of the compact group $Q$ (see~\cite{DR, CDR}). On the other hand, the copy of $\cF(\Aoi)$ in $\cF(\Boi)$ is generated by $\Aoi$ and the Hilbert spaces of isometries implementing the cohomological extensions to $\Boi$ of the covariant DHR endomorphisms with finite statistics of $\Aoi$~\cite[thm.\ 3.5]{CDR}.

Therefore we obtain that $\Boi \subset \cF(\Aoi)$ and then the conclusion follows immediately from~\cite[thm.\ 3.6.a)]{DR} and~\cite[thm.\ 3.4]{CC2}.
\end{proof}

\begin{Remark}\label{rem:Q} According to the general discussion at the beginning of this section, the condition $\Aoi = \cF(\cB)_{0,\iota}^Q$ is automatically satisfied for $Q = G(\cA)_{0,\iota}$ if $\cF(\cA)_{0,\iota} = \cF(\cB)_{0,\iota}$.\end{Remark}

 It would be interesting to know conditions on $\cB$ (and on
$\uo_{0,\iota}$) which guarantee that $\cF(\Boi)$ satisfies the
assumptions in definition~\ref{def:trivialsuper}. It follows from the
discussion in~\cite{DMV} that, since we already assumed that $\Boi$
satisfies Haag duality, assumptions (ii) and (iii) for $\cF(\Boi)$ can
be deduced from analogous assumptions on $\cB$. It is also reasonable to
expect that suitable nuclearity requirements on $\cB$ imply
assumption (iv) for $\cF(\Boi)$, see also the paragraph following proposition~\ref{prop:nosplit}. For what concerns assumption (i), the Reeh-Schlieder property in the scaling limit can be deduced for the algebras $\Boi(W)$ associated to wedges.
Finally, theorem 4.7 of~\cite{CDR} allows to deduce property (v) for $\cF(\Boi)$ from the absence of sectors with infinite statistics for $\Boi$, however it is not clear how to obtain the latter property from the properties of $\cB$.

\subsection{Convergent scaling limits}\label{subsec:gmax}
We now turn to the discussion of the validity of the equality $\Aoi = \cF(\cB)_{0,\iota}^Q$, where actually $Q$ will be a closure of $G(\cA)$ in a suitable topology, and we provide a sufficient condition for it which has some conceptual flavour. In order to do this we appeal to the notion of convergent scaling limit introduced in section~\ref{sec:subsystem}, which is suggested by the experience with models in the perturbative approach to QFT, where there is usually no need of  generalized subsequences in calculating the scaling limit of vacuum expectation values.

We start by showing that, under the additional assumption
that $G(\cB)$ is a normal subgroup of $G(\cA)$, the action of $G(\cA)$ lifts to the scaling algebra $\uFF(\cB)$ and to each scaling limit theory $\cF(\cB)_{0,\iota}$. For simplicity, we also assume that the geometrical symmetry group $\Gamma$ coincides with the translations group, but the arguments below carry over to more general choices.

\begin{Lemma}\label{lem:lift}
Let $\cA \subset \cB$ be an inclusion of Haag dual local nets with $\cH_\cB$ separable and with $G(\cB)$ normal in $G(\cA)$. Then the equation
\begin{equation*}
\ub_\g(\uF)_\l := \b_\g(\uF_\l), \qquad \g \in G(\cA), \,\uF \in \uFF(\cB),\,\l > 0,
\end{equation*}
defines an automorphic action of $G(\cA)$ on $\uFF(\cB)$, which is unitarily implemented in the represention $\poi$ corresponding to any scaling limit state $\ooi$.
\end{Lemma}

\begin{proof}
Let $\uF \in \uFF(\cB)(O)$ and $\g \in G(\cA)$, $g\in G(\cB)$. We get
\begin{equation*}
\sup_{\l > 0}\| \b_g\b_\g(\uF_\l)-\b_\g(\uF_\l)\| = \sup_{\l > 0}\| \b_{\g^{-1}g\g}(\uF_\l)-\uF_\l\|= \|\ub_{\g^{-1}g\g}(\uF)-\uF\|, \\
\end{equation*}
and, since $\g^{-1}g\g \in G(\cB)$ by assumption, the right hand side converges to zero as $g \to e$; analogously
\begin{equation*}
\sup_{\l > 0}\| \a_{\l x}\b_\g(\uF_\l)-\b_\g(\uF_\l)\| = \|\ua_x(\uF)-\uF\|,
\end{equation*}
converges to  $x \to 0$. Therefore we get $\ub_\g(\uF) \in \uFF(\cB)(O)$, and obviously $\g \in G(\cA) \to \ub_\g \in \text{Aut}(\uFF(\cB))$ is a group homomorphism, albeit not pointwise norm continuous. Let then $\uo_{0,\iota}$ be a scaling limit state of $\cF(\cB)$. Since $\o \circ\b_\g = \o$, it follows $\uo_{0,\iota}\circ\ub_\g = \uo_{0,\iota}$, and thus there exists a (not strongly continuous in general) unitary representation $\g \to V_{0,\iota}(\g)$ on $\cH_{\cF(\cB)_{0,\iota}}$ such that $\ad V_{0,\iota}(\g)(\pi_{0,\iota}(\uF))=\pi_{0,\iota}(\ub_\g(\uF))$, where $\pi_{0,\iota}$ is the scaling limit representation associated to $\uo_{0,\iota}$.
\end{proof}

\begin{Theorem}\label{thm:regular}
Let $\cA \subset \cB$ be an inclusion of local nets satisfying the standing assumptions with $G(\cB)$ normal in $G(\cA)$. Moreover, suppose that $\cB$ has a convergent scaling limit and that the algebra $\hat\uBB \subset \uFF(\cB)$ is globally invariant with respect to the action of $G(\cA)$ defined in lemma~\ref{lem:lift}. Then $\Aoi =  \cF(\cB)_{0,\iota}^{G(\cA)}$.
\end{Theorem}

\begin{proof}
The inclusion $\Aoi(O) \subseteq \cF(\cB)_{0,\iota}(O)^{G(\cA)}$ is trivial: given $A = \pi_{0,\iota}(\uA)$, $\uA \in \uAA(O)$, we obviously have $\uA_\l \in \cF(\cB)(\l O)^{G(\cA)}$ and therefore, by definition of $\ub$, $\ub_\g(\uA) = \uA$, which entails $A \in \cF(\cB)_{0,\iota}(O)^{G(\cA)}$.

In order to prove the converse inclusion, let $F \in \cF(\cB)_{0,\iota}(O)^{G(\cA)}$, and note that, since $ \cF(\cB)_{0,\iota}(O)^{G(\cA)} \subset \Boi(O)$, we can choose elements $\uF_n \in \hat\uBB$ such that $\poi(\uF_n)$ converges strongly to $F$ as $n \to +\infty$. We define
\begin{equation*}
\uA_{n,\l} := \int_{G(\cA)} d\g\,\g\uF_{n,\l}\g^{-1},
\end{equation*}
where the integral is defined in the weak topology. It is plain that $\uA_{n,\l} \in \cA(\l O)$ and $\| \uA_{n,\l} \| \leq \|\uF_n\|$, and furthermore
\begin{equation*}
\| \a_{\l x}(\uA_{n,\l})-\uA_{n,\l} \| = \left\| \int_{G(\cA)} d\g\,\g\left(\a_{\l x}(\uF_{n,\l})-\uF_{n,\l}\right)\g^{-1}\right\| \leq \| \ua_x(\uF_n)-\uF_n\|,
\end{equation*}
which gives $\uA_n \in \uAA(O)$. Moreover, for fixed $n \in \bN$, it is possible to find a sequence $(\l_m)_{m\in\bN}$ such that
\begin{equation}\label{eq:FmenoAseq}\begin{split}
\| [\pi_{0,\iota}(\uF_n)-\pi_{0,\iota}(\uA_n)]\O_{0,\iota}\| &= \lim_{m\to+\infty} \left\| \int_{G(\cA)} d\g\,\left(\uF_{n,\l_m}-\g \uF_{n,\l_m}\right)\O\right\| \\
&\leq \lim_{m \to +\infty} \int_{G(\cA)} d\g\,\left\|\left(\uF_{n,\l_m}-\g \uF_{n,\l_m}\right)\O\right\| .
\end{split}\end{equation}
Thanks to the $G(\cA)$-invariance of $\hat\uBB$, $\ub_\g (\uF_n)\in\hat\uBB$ so that, for each $\gamma \in G(\cA)$, there holds
\begin{equation}\label{eq:FmenogammaFlambda}\begin{split}
\lim_{\l \to 0} \left\|\left(\uF_{n,\l}-\g \uF_{n,\l}\right)\O\right\|^2 &= \lim_{\l\to 0}\o\Big(\left(\uF_{n,\l}-\ub_\g( \uF_n)_\l\right)^*\left(\uF_{n,\l}-\ub_\g( \uF_n)_\l\right)\Big) \\
&= \big\|\big[\pi_{0,\iota}(\uF_n)-\ad V_{0,\iota}(\g)(\pi_{0,\iota}(\uF_n))\big]\O_{0,\iota}\big\|^2\\
&\leq 4 \| [\pi_{0,\iota}(\uF_n)-F]\O\|^2 .
\end{split}\end{equation}
Therefore, by applying Lebesgue's dominated convergence theorem to~\eqref{eq:FmenoAseq}, we conclude that
\begin{equation*}
\| [\pi_{0,\iota}(\uF_n)-\pi_{0,\iota}(\uA_n)]\O_{0,\iota}\|\leq 2 \| [\pi_{0,\iota}(\uF_n)-F]\O\|,
\end{equation*}
which, together with the fact that $\O_{0,\iota}$ is separating for the local algebras, gives us that $\pi_{0,\iota}(\uA_n)$ converges strongly to $F$ as $n \to +\infty$.
\end{proof}

\begin{Corollary}\label{cor:reg_fixpoint}
Under the same assumptions as in theorem~\ref{thm:regular}, and assuming that $\cF(\Boi)$ satisfies the properties (i)-(v) in definition \ref{def:trivialsuper}, there holds $\cF(\Aoi) = \cF(\Boi)$, and $\Aoi$ is the fixpoint net of $\cF(\Boi)$ under a compact group of internal symmetries.
\end{Corollary}

\begin{proof}
If $Q$ denotes the closure, in the strong operator topology on  $\cH_{\cF(\cB)_{0,\iota}}$, of the group of unitaries $\{\Voi(\g)\,:\,\g \in G(\cA)\}$, it is an easy consequence of theorem~\ref{thm:regular} that $\Aoi = \cF(\cB)_{0,\iota}^Q$. Furthermore $\cF(\cB)_{0,\iota} \subset \cF(\cB_{0,\iota})$ satisfies the split property~\cite[sec.\ 2]{CC2}, and therefore $Q$ is compact~\cite{DL}, so that theorem~\ref{thm:Q} gives the statement.
\end{proof}

The fact that the scaling limit is convergent has been checked in~\cite{BDM} for the theory of a single massive free scalar field. The $G(\cA)$-invariance condition used in the above theorem can possibly be shown for the theory of a multiplet of free scalar fields in the following way. If $\cF$ is the field net generated by such a multiplet, it should be possible,
using the same techniques as in~\cite{BDM},
to construct a C$^*$-subalgebra $\hat \uFF \subset \uFF$ which is globally invariant under $G(\cA) = G_\text{max}$, the maximal group of internal symmetries of $\cF$,\footnote{$G_\text{max}$ is the group of unitaries $U$ on $\cH_\cF$ such that $U \cF(O) U^* = \cF(O)$, $U U(\g) = U(\g) U$ for each $\g \in \Gamma$, and $U \O = \O$.} and on which a given normal subgroup $G(\cB) \subset G(\cA)$ acts strongly continuously. Then if $\hat \uBB := \smash{\hat \uFF}^{G(\cB)} \subset \uBB$, one has that $\hat \uBB$ is $G(\cA)$-invariant thanks to the normality of $G(\cB)$ in $G(\cA)$. Moreover, $\hat\uBB$ has the two properties required in the definition of a convergent scaling limit. It is plain that there exists $\lim_{\l\to 0}\o(\uB_\l)$ for each $\uB \in \hat\uBB$, as this holds for $\hat\uFF$, while the property $\poi(\hat\uBB(O))'' = \Boi(O)$ follows from the analogous property for $\hat\uFF$ by averaging in the usual way with respect to the strongly continuous action of $G(\cB)$.

Finally, we notice that the fact that the scaling limit is convergent for the free scalar field depends on the nuclearity properties of the theory.

\bigskip
Without the assumption of a convergent scaling limit, the above proof breaks down because of the necessity of interchanging the integral on $G(\cA)$ in equation~\eqref{eq:FmenoAseq} with the limit along a generalized sequence $(\l_\k)_\k$, which is not  guaranteed under the present conditions. One can only speculate that additional assumptions (e.g.\ nuclearity) may provide further insight on this issue. Anyway, if one cannot take the limit under the integral sign in the above discussion, we have to leave open the possibility that $\Aoi \subsetneq \cF(\cB)_{0,\iota}^{G(\cA)}$, in which case we are left with two mutually exclusive possibilities: either there exists some compact group $Q$ ``larger'' than
 (the strong operator closure of)
 $G(\cA)$ acting on $\cF(\cB)_{0,\iota}$ such that $\Aoi = \cF(\cB)_{0,\iota}^Q$, or there is no such group. In the former case the principle of gauge invariance is restored at the price of ``enlarging'' the group $G(\cA)$. As an illustration of the physical meaning of such situation, consider the particular case in which $\cA = \cF^G$ (i.e.\ $G(\cB)$ is trivial and $G=G(\cA)$, $\cF = \cF(\cA)=\cF(\cB)=\cB$): then the existence of $Q$ would mean that $\Aoi$ is the fixed point net of the ``wrong'' scaling limit field net $\Foi$, defined without any reference to the action of $G$, and
 would imply that it is possible to create ``new''
sectors of $\Aoi$ by looking at the scaling limit of states where a region of radius $\l$ contains an amount of charge which grows unboundedly as $\l \to 0$. Such sectors should however not be regarded as confined, as they could be created by performing suitable operations at finite, albeit small, scales. A thorough analysis of the structure of such sectors is of considerable interest in itself, and would require going beyond the framework of~\cite{DMV}.

\subsection{On the scaling of Noether currents}
As an application of the above results, we discuss the scaling limit of nets generated by local implementations of symmetries~\cite{BDoL}. Let $\cB$ be a local net and let $\cA\subset\cB$ be the dual of the net generated by the local implementations of translations of $\cF(\cB)$. The validity of the equality $\cA = \cF(\cB)^{G_\text{max}}$, where $G_\text{max}$ is the maximal group of internal symmetries of $\cF(\cB)$, has been thoroughly discussed in~\cite{CC1,CC2}.

\begin{Theorem}\label{thm:noether}
Let $\cA \subset \cB$  satisfy the standing assumptions, where $\cA$ is the dual of the net generated by the canonical local implementers of the translations of $\cF(\cB)$, and suppose furthermore that $\cF(\Boi)$ satisfies the properties (i)-(v) in definition \ref{def:trivialsuper}, that $G(\cB)$ is normal in $G_\textup{max}$, and that $\cB$ has a convergent scaling limit such that $\hat\uBB$ is $G_\textup{max}$-invariant. Let $\tilde \cA_{0,\iota}$ be the dual of the net generated by the local implementations of translations of $\cF(\Boi)$. Then
$$\tilde\cA_{0,\iota}\subset \Aoi.$$
\end{Theorem}

\begin{proof}
Thanks to corollary~\ref{cor:reg_fixpoint}, one has $\cF(\Aoi) = \cF(\Boi)$, and $\Aoi = \cF(\Boi)^{G(\Aoi)}$, and therefore, with $\tilde G_\text{max}$ the maximal group of internal symmetries of $\cF(\Boi)$, $\tilde \cA_{0,\iota} = \cF(\Boi)^{\tilde G_\text{max}} \subset \Aoi$.
\end{proof}

In short, the above result states that the scaling limit of the net generated by the local energy-momentum tensor contains the net generated by the local energy-momentum tensor of the scaling limit.\footnote{This conclusion is supported by some preliminary calculations performed directly on the universal localizing maps that are used to construct the canonical local implementers.} It is likely that in favourable circumstances $\tilde\cA_{0,\iota} = \Aoi$. This is trivially illustrated by the example of the fixpoints of the free field net discussed in the previous section.

In the case in which the scaling limit of $\cB$ is not convergent, one has to look back at theorem~\ref{thm:Q}.
In turn, one should be able to show by similar methods an analogous result for more general Noether currents,
cf.\cite{CC2}.

\medskip
An issue that should be taken into account is the fact that the split property is not necessarily preserved
in the scaling limit. This somehow unpleasant feature, although strictly speaking ruled out by our assumptions,
can partly justify at a heuristic level
the possibly strict inclusion of nets
that we obtained in theorem~\ref{thm:noether}.
In fact, in that case one cannot even define the local implementers of the scaling limit although it still makes sense
to consider the scaling limit of the net generated by the Noether charges of the original theory.

Examples of local nets satisfying the split property but whose scaling limit does not satisfy it can be easily found.

\begin{Proposition}\label{prop:nosplit}
Let $\cB$ be the dual of the local net generated by a generalized free field with a mass measure $d\rho(m) = \sum_i \d(m-m_i)$, such that
\begin{equation}\label{eq:sum}
\sum_i  e^{-\gamma m_i} < \infty
\end{equation}
for all $\gamma>0$. Then the split property holds for $\cB$ but for none of its scaling limit nets $\Boi$.
\end{Proposition}

\begin{proof}
Since condition~\eqref{eq:sum} implies $\sum_i m_i^4 e^{-\d m_i}<\infty$ for each $\d > 0$, the split property for $\cB$ follows from~\cite[p.\ 529]{DL} and~\cite[cor.\ 4.2]{DDFL}. By~\cite[thm.\ 4.1]{Moh} each scaling limit net $\Boi$ (contains a subnet that) does not satisfy the Haag-Swieca compactness condition, and thus, by~\cite[prop.\ 4.2]{BDL}, it does not satisfy the split property either.
\end{proof}

A related problem is to provide conditions on $\cB$ ensuring that $\cF(\Boi)$ enjoys the split property. By the results in~\cite[cor.\ 4.6]{BD}, some conditions on $\Boi$ are known to imply suitable nuclearity properties of $\cF(\Boi)$, which in turn imply the split property by~\cite[sec.\ 4]{BDL}. On the other hand, the methods employed in~\cite[thm.\ 4.5]{Buc} to prove nuclearity properties of the scaling limit theory $\Boi$ starting from certain phase space behaviour of the underlying theory $\cB$ can possibly be adapted to show that (some of) the conditions on $\Boi$ considered in~\cite{BD} follow from appropriate nuclearity requirements on $\cB$.

\subsection{Preserved sectors}
The notion of preserved DHR sector has been introduced in~\cite[def.\ 5.4]{DMV}. In the spirit of the present paper a natural question concerns the relation between the preservation of DHR morphisms with finite statistics of $\cA$ and $\cB$. Clearly, the cohomological extension property of morphisms plays again a crucial role. 
In particular, it is reasonable to expect
that if all the morphisms of $\cB$ are preserved, then the same will be true for the morphisms of $\cA$, since the Hilbert spaces of isometries in $\cF(\cB)$ implementing the cohomological extension to $\cB$ of a given morphism of $\cA$ would satisfy the preservation condition.
Possible applications of such result include a generalization of the theorem in the previous section, where we replace $\cF(\cA)_{0,\iota}$ and $\cF(\cB)_{0,\iota}$ with the subnets generated by the isometries associated to the (scaling limits of the) preserved sectors, which should be automatically independent of the gauge groups, and therefore coincide.

However, what is missing in the above argument is the fact that the analysis in \cite{DMV}
has been carried out only for irreducible morphisms (while the extension maps in general irreducible morphisms to reducible ones) and,
although there is no apparent obstruction for extending it to the reducible case,
in the reminder of this short section
we will limit ourselves to some simple remarks, postponing a thorough analysis of this point to future work.

\medskip
We consider the situation outlined at the beginning of section~\ref{subsec:general}, i.e.\ an inclusion $\cA \subset \cB$ of Haag dual and Poincar\'e covariant nets in the vacuum sector with $\cF(\cA) = \cF(\cB)$, and scaling limit states $\uooi^\cB$ of $\cB$ and $\uooi^\cA = \uooi^\cB \rest \uAA$ of $\cA$. In the following result we make use of the notion of asymptotic containment, introduced in~\cite[def.\ 5.2]{DMV}.

\begin{Proposition}
Let $\xi$ be a $\uooi^\cA$-preserved class of DHR morphisms of $\cA$,
and let  $\psi_j(\l) \in \cF(\l O)$ be an associated scaled multiplet which is asymptotically contained in $\uFF(\cA)$. Then the cohomological extension $\hat{\xi}$ of $\xi$ to $\cB$ is
$\uooi^\cB$-preserved, with $\psi_j(\l)$ an associated scaled multiplet asymptotically contained in $\uFF(\cB)$.
\end{Proposition}

\begin{proof}
If $\rho_\l$ is the DHR morphism of $\cA$ in the class $\xi$ localized in $\l O$ implemented by the multiplet $\psi_j(\l)$, then its
cohomological extension $\hat{\rho}_\l$ is in the class $\hat{\xi}$, still localized in $\l O$ and also implemented by $\psi_j(\l)$.
Now, since $\uFF(\cA) \subset \uFF(\cB)$, it is immediate to conclude that $\psi_j(\l)$ is asymptotically contained in $\uFF(\cB)$, and therefore $\hat{\xi}$ is $\uooi^\cB$-preserved.
\end{proof}

We define $\cF(\cA)^{\text{pres}}_{0,\iota}$ as the net generated
by $\Aoi$ and the scaling limits of scaled multiplets asymptotically contained in $\uFF(\cA)$ associated to $\uooi^\cA$-preserved sectors of $\cA$, see prop.\ 5.5 in~\cite{DMV}.  Likewise we define $\cF(\cB)^{\text{pres}}_{0,\iota}$ with respect to the scaling limit state $\uooi^\cB$.

\begin{Corollary}
With the above notations, 
there holds
\begin{equation*}
\begin{array}{ccccc}
\Aoi &\subset &\cF(\cA)^{\textup{pres}}_{0,\iota}          &\subset &\cF(\cA)_{0,\iota} \\
\cap          &        &\cap                                 &        &\cap \\
\Boi &\subset &\cF(\cB)^{\textup{pres}}_{0,\iota}  &\subset &\cF(\cB)_{0,\iota}
\end{array}
\end{equation*}
\end{Corollary}

In some cases one has $\cF(\cA)^{\text{pres}}_{0,\iota} = \cF(\cB)^{\text{pres}}_{0,\iota}$. For free fields this follows from example~\ref{ex:free} and the fact that all sectors of the fixpoint net of the free field are preserved.
Another example is given by a dilation invariant theory satisfying the Haag-Swieca compactness condition, where it follows
easily from the results of~\cite{BV, DMV}, that $\cF(\cA)^{\text{pres}}_{0,\iota} = \cF(\cB)^{\text{pres}}_{0,\iota} = \cF(\cA)_{0,\iota} = \cF(\cB)_{0,\iota} = \cF$.

\section{Final comments}
We end this paper with few comments on further possible extension of the results presented above, in addition to those already mentioned in the main body.

\medskip
Given a subsystem $\cA \subset \cB$ as in section~\ref{sec:DR}, we assumed that $\cF(\cA) = \cF(\cB)$. However in general it holds $\cF(\cA) \subset \cF(\cB)$~\cite{CDR} and, in the situation considered in~\cite{CC1, CC2}, it is shown that $\cF(\cA) = \cF(\cB)\otimes \cC$ (graded tensor product)
for a suitable net $\cC$ . Therefore, in order to treat this more general framework, one should generalize the results about the scaling limit of tensor product theories in~\cite{DM}.

\medskip
Another natural example of subsystem, to which most of our results don't apply, is provided by the inclusion $\cA \subset \cA^d$ of a net into its dual,
a situation that arises typically when there are spontaneously broken symmetries. The analysis of
the structure of such subsystems in the scaling limit has some interest as it could possibly simplify the
study of the relations between the superselection structures of $\cA$ and of $\Aoi$. For instance, sufficient conditions on $\cA$ which imply essential duality, but not duality, of $\Aoi$ are known, so it would be interesting to know when the scaling limit of $\cA^d$ coincides with the dual of $\Aoi$.

\medskip
We conclude by mentioning few very intriguing but rather speculative ideas.
In \cite{DMV} it has been shown that it is possible to formulate conditions on the
scaling limit of a theory which imply the equality of local
and global intertwiners.
There are other long-standing structural problems in superselection theory that could hopefully be related one way or another
to the short distance properties of the theory. Just to give some example, we cite here the problem of
recovering pointlike
Wightman fields with specific physical interpretation out of local algebras (i.e., a full quantum Noether theorem), and that of
ruling out the existence of sectors
with infinite statistics.


\begin{thebibliography}{99}

\bibitem{BDM1}H. Bostelmann, C. D'Antoni, G. Morsella, Scaling algebras and pointlike fields. A nonperturbative approach to renormalization, arXiv:0711.4237, to appear on Comm. Math. Phys.
\bibitem{BDM}H. Bostelmann, C. D'Antoni, G. Morsella, work in progress.
\bibitem{Buc1}D. Buchholz, Quarks, gluons, colour: facts or fiction? Nucl. Phys. B 469 (1996) no. 1-2, 333-353.
\bibitem{Buc}D. Buchholz, Phase space properties of local observables and structure of scaling limits. Ann. Inst. H. Poincar\'e Phys. Theor.  64  (1996),  no. 4, 433--459.
\bibitem{BD}D. Buchholz, C. D'Antoni, Phase space properties of charged fields in theories of local observables. Rev. Math. Phys.  7  (1995),  no. 4, 527--557.
\bibitem{BDL}D. Buchholz, C. D'Antoni, R. Longo, Nuclear maps and modular structures. II. Applications to quantum field theory.  Comm. Math. Phys.  129  (1990),  no. 1, 115--138.
\bibitem{BDoL}D. Buchholz, S. Doplicher, R. Longo, On Noether's theorem in quantum field theory. Ann. Phys. 170 (1986), 1--17.
\bibitem{BV} D. Buchholz, R. Verch,
Scaling algebras and renormalization group in algebraic quantum field theory.
Rev. Math. Phys.  7  (1995),  no. 8, 1195--1239.
\bibitem{BV2} D. Buchholz, R. Verch,
Scaling algebras and renormalization group in algebraic quantum field theory. II.
Instructive examples.  Rev. Math. Phys.  10  (1998),  no. 6, 775--800.
\bibitem{CC1} S. Carpi, R. Conti,
Classification of subsystems for local nets with trivial superselection structure.
Comm. Math. Phys.  217  (2001),  no. 1, 89--106.
\bibitem{CC2} S. Carpi, R. Conti,
Classification of subsystems for graded-local nets with trivial superselection
structure.  Comm. Math. Phys.  253  (2005),  no. 2, 423--449.
\bibitem{CDR}R. Conti, S. Doplicher, J. E. Roberts,
Superselection theory for subsystems.  Comm. Math. Phys.  218  (2001),  no. 2,
263--281.
\bibitem{DDFL}C. D'Antoni, S. Doplicher, K. Fredenhagen, R. Longo, Convergence
 of local charges and continuity properties of $W\sp *$-inclusions.  Comm. Math. Phys.  110  (1987),  no. 2, 325--348.
\bibitem{DMV}C. D'Antoni, G. Morsella, R. Verch,
Scaling algebras for charged fields and short-distance analysis for localizable
and topological charges.  Ann. Henri Poincar\'e  5  (2004),  no. 5, 809--870.
\bibitem{DM}C. D'Antoni, G. Morsella,
Scaling algebras and superselection sectors: study of a class of models.  Rev.
Math. Phys.  18  (2006),  no. 5, 565--594.
\bibitem{DL}S. Doplicher, R. Longo, Standard and split inclusions of von Neumann algebras.
Invent. Math.  75  (1984),  no. 3, 493--536.
\bibitem{DR}S. Doplicher, J. E. Roberts,
Why there is a field algebra with a compact gauge group describing the
superselection structure in particle physics.  Comm. Math. Phys.  131  (1990),
no. 1, 51--107.
\bibitem{Dri}W. Driessler, Duality and absence of locally generated superselection sectors for CCR-type algebras, Comm. Math. Phys. 70 (1979), no. 3, 213-220.
\bibitem{Haag}R. Haag, \emph{Local quantum physics}, IInd edition, Springer-Verlag, Berlin, 1996.
\bibitem{LR}R. Longo, K.-H. Rehren, Nets of subfactors, Rev. Math. Phys. 7 (1995) no. 4, 567-597.
\bibitem{L}M. Lutz, Ein lokales Netz ohne Ultraviolettfixpunkte der Renormierungsgruppe,
diploma thesis, Hamburg University (1997).
\bibitem{Moh} S. Mohrdieck,
Phase space structure and short distance behaviour of local quantum field theories.
J. Math. Phys. 43 (2002), 3565-3574.
\bibitem{Rob}J. E. Roberts, Localization in algebraic field theory, Comm. Math. Phys.
85 (1982), 87-98.
\bibitem{Wi}E. H. Wichmann, On systems of local generators and the duality condition, J. Math Phys. 24 (1983), 1633-1644.
\bibitem{Y}K. Yosida, \emph{Functional analysis},
Reprint of the sixth (1980) edition. Classics in Mathematics.
Springer-Verlag, Berlin, 1995.
\end{thebibliography}
\end{document}